\documentclass[a4paper,11pt,reqno]{amsart}

\usepackage{dsfont,amssymb}

\numberwithin{equation}{section}

\newtheorem{theorem}{Theorem}[section]

\newtheorem{lemma}[theorem]{Lemma}
\newtheorem{proposition}[theorem]{Proposition}
\theoremstyle{remark}

\newtheorem{example}[theorem]{Example}
\theoremstyle{definition}
\newtheorem{definition}[theorem]{Definition}

\newcommand\bp{\begin{proof}}
\newcommand\ep{\end{proof}}
\newcommand\ee{\nopagebreak\mbox{\ }\hfill$\diamondsuit$}

\newcommand\Dhat{{\hat\Delta}}

\newcommand\Ad{\operatorname{Ad}}

\newcommand{\Ghop}{{\hat G}^{\rm op}}
\newcommand{\Ghopo}{({\hat G_\Omega})^{\rm op}}
\newcommand{\Gop}{{G}^{\rm op}}
\newcommand{\Dhop}{{\hat\Delta}^{\rm op}}
\newcommand{\Dop}{{\Delta}^{\rm op}}
\newcommand{\Whop}{{\hat W}^{\rm op}}

\newcommand{\C}{{\mathbb C}}
\newcommand{\R}{{\mathbb R}}
\newcommand\T{{\mathbb T}}
\newcommand\Z{{\mathbb Z}}

\newcommand{\A}{{\mathcal A}}
\newcommand{\E}{{\mathcal E}}
\renewcommand{\H}{{\mathcal H}}
\newcommand{\F}{{\mathcal F}}

\newcommand{\PP}{{\mathcal P}}

\newcommand\g{{\mathfrak g}}

\newcommand\enu[1]{\smallskip\newline\makebox[5mm][l]{\rm(#1)}}

\begin{document}

\title[Deformation of C$^*$-algebras]{Deformation of C$^*$-algebras by cocycles on locally compact quantum groups}

\author[S. Neshveyev]{Sergey Neshveyev}

\email{sergeyn@math.uio.no}

\address{Department of Mathematics, University of Oslo,
P.O. Box 1053 Blindern, NO-0316 Oslo, Norway}

\author[L. Tuset]{Lars Tuset}

\email{Lars.Tuset@hioa.no}

\address{Department of Computer Science, Oslo and Akershus University College of Applied Sciences,
P.O. Box 4 St. Olavs plass, NO-0130 Oslo, Norway}

\thanks{The research leading to these results has received funding from the Research Council of Norway and the European Research Council
under the European Union's Seventh Framework Programme (FP/2007-2013) / ERC Grant Agreement no. 307663
}

\date{January 21, 2013; revised version December 21, 2013}

\begin{abstract}
Given a C$^*$-algebra $A$ with a left action of a locally compact quantum group $G$ on it and a unitary $2$-cocycle $\Omega$ on $\hat G$, we define a deformation $A_\Omega$ of $A$. The construction behaves well under certain additional technical assumptions on $\Omega$, the most important of which is regularity, meaning that $C_0(G)_\Omega\rtimes G$ is isomorphic to the algebra of compact operators on some Hilbert space. In particular, then $A_\Omega$ is stably isomorphic to the iterated twisted crossed product $\Ghop\ltimes_\Omega G\ltimes A$. Also, in good situations, the C$^*$-algebra $A_\Omega$ carries a left action of the deformed quantum group~$G_\Omega$ and we have an isomorphism $G_\Omega\ltimes A_\Omega\cong G\ltimes A$. When $G$ is a genuine locally compact group, we show that the action of $G$ on $C_0(G)_\Omega=C^*_r(\hat G;\Omega)$ is always integrable. Stronger assumptions of properness and saturation of the action imply regularity. As an example, we make a preliminary analysis of the cocycles on the duals of some solvable Lie groups recently constructed by Bieliavsky et~al., and discuss the relation of our construction to that of Bieliavsky and Gayral.
\end{abstract}

\maketitle

\section*{Introduction}

Assume $(\H,\Dhat)$ is a Hopf algebra and $\A$ is a left $\H$-module algebra, with the action of $\H$ denoted by $h\otimes a\mapsto h\rhd a$. Assume also that $\Omega\in\H\otimes\H$ is an invertible element satisfying the cocycle identity
$$
(\Omega\otimes1)(\Dhat\otimes\iota)(\Omega)=(1\otimes\Omega)(\iota\otimes\Dhat)(\Omega).
$$
In this case we can consider a new Hopf algebra $\H_\Omega$, defined by Drinfeld~\cite{Dr}, such that $\H_\Omega=\H$ as an algebra, but the coproduct is given by $\Dhat_\Omega=\Omega\Dhat(\cdot)\Omega^{-1}$. We can also define a new product~$\star_\Omega$ on~$\A_\Omega=\A$~by
$$
a\star_\Omega b=m(\Omega^{-1}\rhd(a\otimes b)),
$$
where $m(a\otimes b)=ab$. Then $\A_\Omega$ is an $\H_\Omega$-module algebra and, as was observed by Majid~\cite{Ma} (see also~\cite{BPVO}), for the corresponding smash, or crossed, products we have
\begin{equation} \label{eMa}
\H_\Omega\# \A_\Omega\cong\H\#\A.
\end{equation}

Our goal in this paper is to develop a similar theory in the context of C$^*$-algebras and actions of locally compact quantum groups. Thus, given a C$^*$-algebra $A$ with an action of a locally compact quantum group $G$ on it, and a unitary cocycle $\Omega\in L^\infty(\hat G)\bar\otimes L^\infty(\hat G)$, we want to define a deformation~$A_\Omega$ of~$A$. Note that the deformed quantum group $G_\Omega$ is defined in full generality by the theory developed by De Commer~\cite{DC}.

Particular cases of our construction of $A_\Omega$ are of course well-known. For example, when $G$ is dual to a discrete group $\Gamma$, $A=\Gamma\ltimes_{\gamma} B$ (all crossed products in this paper are assumed to be reduced) and the action of $G$ is the dual action $\hat \gamma$, then $A_\Omega$ is nothing else than the twisted crossed product $\Gamma\ltimes_{\gamma,\Omega}B$, as defined already by Zeller-Meier~\cite{ZM} in the 60s.

In the case when $G$ is a compact group, a study of cocycles on $\hat G$ was initiated by Landstad~\cite{La} and Wassermann~\cite{WaPhD,Wa} in the early 80s. They defined twisted group C$^*$-algebras $C^*_r(\hat G;\Omega)$, which should be thought of as deformations $C_0(G)_\Omega$.

Another milestone is the work of Rieffel~\cite{Ri1} for $G=\R^{2d}$. He was motivated by deformation quantization theory and extended Weyl quantization to actions of $\R^{2d}$ on C$^*$-algebras. His theory is beautiful, but quite complicated, based on an extension of oscillatory integrals to C$^*$-algebras. A much simpler, although less explicit, approach was later proposed by Kasprzak~\cite{Kas}. His idea was to start with isomorphism~\eqref{eMa}. This isomorphism implies that $\A_\Omega$ can be identified with the fixed point subalgebra of $\H\#\A$ with respect to a coaction of $\H_\Omega$ that corresponds to the dual coaction on $\H_\Omega\#\A_\Omega$. It is easy to check that this coaction of $\H_\Omega$ on $\H\#\A$ is simply the dual coaction of $\H$ twisted by $\Omega$. This allows one to describe $\A_\Omega$ in terms of $\H\#\A$ and $\Omega$ without using expressions like $m(\Omega^{-1}\rhd(a\otimes b))$ that are difficult to makes sense of in the analytic setting. Kasprzak developed this idea in the  setting of C$^*$-algebras in the case when $G$ is an abelian locally compact group and~$\Omega$~is a continuous $2$-cocycle on $\hat G$. But his theory works equally well when $G$ is the dual of a locally compact group.

In the previous paper the first author together with Bhowmick and Sangha~\cite{BNS} extended Kasprzak's approach to the case of measurable cocycles. The problem with such cocycles is that the twisted dual action of $\hat G$ on $G\ltimes A$ (since we deal now with group duals the deformed quantum group $G_\Omega$ is $G$) is not always well-defined. Nevertheless, a description of $A_\Omega$ in terms of generators continues to make sense and, as was shown in~\cite{BNS}, these algebras still satisfy a number of properties to be considered as the correct deformations of~$A$.

A by no means exhaustive list of other relevant papers on cocycle deformations in the operator algebraic setting includes \cite{BDRV}, \cite{Betal}, \cite{BG}, \cite{CL}, \cite{DC2}, \cite{FV}, \cite{HM}, \cite{Kas2}, \cite{LR1}, \cite{LR2}, \cite{Si}, \cite{VV}, \cite{Va}, \cite{Wa1}, \cite{Wa2}, \cite{Y}.

\smallskip

In this paper we continue the work started in \cite{BNS} and define the deformations $A_\Omega$ for general locally compact quantum groups $G$ and arbitrary unitary cocycles $\Omega$ on $\hat G$. In fact, our primary interest is the group case. This is the situation studied by Landstad and Raeburn~\cite{LR1,LR2} for C$^*$-algebras of the form $A=C_0(G/H)$ and a particular class of cocycles on $\hat G$. But since, as follows from the above discussion, the deformed quantum group $G_\Omega$ should play a role in the theory, trying to work only with groups and their duals is unnecessarily restrictive. Moreover, in the proofs of a significant number of general results there would be almost no simplifications even if we restricted ourselves to the group situation.

The paper is organized as follows. In Section~\ref{s1} we collect some basic facts on locally compact quantum groups.

In Section~\ref{s2} we study various notions related to cocycles. In particular, here we introduce twisted crossed products. They are related, but not in the most straightforward way, to cocycle crossed products studied by Vaes and Vainerman~\cite{VV}. Another important notion is regularity of a cocycle, which means that the crossed product $C^*_r(\hat G;\Omega)\rtimes G$ is isomorphic to the algebra of compact operators on $L^2(G)$. For regular quantum groups we show that regularity of a cocycle $\Omega$ is equivalent to the inclusion $(K\otimes1)\hat W\Omega^*(1\otimes K)\subset K\otimes K$, where $K$ is the algebra of compact operators on $L^2(G)$ and $\hat W$ is the multiplicative unitary of~$\hat G$. But we leave open the question of finding somewhat more manageable sufficient conditions for regularity.

Section~\ref{s3} contains our main general results. Here we introduce the deformed algebras $A_\Omega$ and study such questions as the relation of $A_\Omega$ to twisted crossed products, existence of an action of~$G_\Omega$ on~$A_\Omega$, deformation in stages, invariance of $A_\Omega$ under replacing $\Omega$ by a cohomologous cocycle.

In Section~\ref{sGD} we specialize to the group case. The main goal is to understand when a cocycle on~$\hat G$ is regular, but the outcome is far from satisfactory. We observe that $\Omega$ is regular if the action of~$G$ on $C^*_r(\hat G;\Omega)$ is proper and saturated in the sense of Rieffel~\cite{Ri}. What we are able to prove in general is that this action always has a weaker property of integrability; note that integrability in an even weaker sense has already been established by Vaes and Vainerman~\cite{VV}. One outcome of this discussion is that $A_\Omega$ is generated by the image of a dense subspace of $C^*_r(\hat G;\Omega)\otimes A$ under the operator-valued weight $C^*_r(\hat G;\Omega)\otimes A\to M(C^*_r(\hat G;\Omega)\otimes A)^G$.

In Section~\ref{sBG}, in order to illustrate some of our general results, as well as the difficulties one might encounter in analyzing concrete examples, we briefly consider the cocycles on the duals of some solvable Lie groups recently constructed by Bieliavsky et~al.~\cite{BG,Betal}. A detailed study will appear elsewhere.

We finish the paper with a list of open problems in Section~\ref{sOP}.

\smallskip

\noindent{\bf Acknowledgement.} We would like to thank Jyotishman Bhowmick for fruitful discussions and a careful reading of the manuscript. We are also grateful to Kenny De Commer for comments on the first version of the paper.

\bigskip

\section{Preliminaries} \label{s1}

In this section we will fix our notation and recall some facts on locally compact quantum groups that we will use freely throughout the paper.

\subsection{Locally compact quantum groups}

Recall~\cite{KV,KVvN} that a locally compact quantum group, in the von Neumann algebraic setting, is a pair $G=(M,\Delta)$ consisting of a von Neumann algebra~$M$ and a coassociative normal unital $*$-homomorphism $\Delta\colon M\to M\bar\otimes M$ such that there exist a left invariant n.s.f weight $\varphi$ and a right invariant n.s.f. weight $\psi$ on $M$. We will often use the suggestive notation $L^\infty(G)$ for $M$. Denote by $L^2(G)$ the space of the GNS-representation of $M$ defined by the left invariant Haar weight $\varphi$. Write $\Lambda\colon \mathcal N_\varphi\to L^2(G)$ for the corresponding map, where $\mathcal N_\varphi=\{x\in L^\infty(G)\mid\varphi(x^*x)<\infty\}$. Then the multiplicative unitary $W$ of $G$ is defined by
$$
W^*(\Lambda(x)\otimes\Lambda(y))=(\Lambda\otimes\Lambda)(\Delta(y)(x\otimes1)),\ \ \text{for}\ \ x,y\in\mathcal N_\varphi.
$$
Therefore, identifying $L^\infty(G)$ with its image under the GNS-representation defined by $\varphi$, we have
$$
\Delta(x)=W^*(1\otimes x)W\ \ \text{for}\ \ x\in L^\infty(G).
$$

Throughout the whole paper we will denote by $K$ the C$^*$-algebra of compact operators on $L^2(G)$. We identify $K^*$ with $B(L^2(G))_*$. For a subset $X$ of a normed space we denote by $[X]$ the norm closure of the linear span of $X$. Using this notation the C$^*$-algebra $C_0(G)$ of continuous functions on $G$ vanishing at infinity is defined by
$$
C_0(G)=[(\iota\otimes\omega)(W)\mid \omega\in K^*].
$$

The dual quantum group $\hat G=(\hat M,\Dhat)$ is defined by
$$
\hat M=\{(\omega\otimes\iota)(W)\mid\omega\in K^*\}^{-\sigma\text{-strong}^*},\ \ \Dhat(x)=\Sigma W(x\otimes1)W^*\Sigma,
$$
where $\Sigma$ is the flip on $L^2(G)\otimes L^2(G)$. By definition $\hat M$ is represented on $L^2(G)$. This representation is identified with the GNS-representation defined by a left invariant Haar weight $\hat\varphi$ on $\hat G$, with the corresponding map $\hat\Lambda\colon\mathcal N_{\hat\varphi}\to L^2(G)$ uniquely defined by the identities
$$
(\hat\Lambda((\omega\otimes\iota)(W)),\Lambda(x))=\omega(x^*)
$$
for $x\in\mathcal N_\varphi$ and suitable $\omega\in K^*$, namely, for $\omega$ such that the map $\Lambda(x)\mapsto\overline{\omega(x^*)}$ extends to a bounded linear functional on $L^2(G)$. Under this identification the multiplicative unitary $\hat W$ of $\hat G$ is given by
$$
\hat W=\Sigma W^*\Sigma.
$$
We then have
$$
C_0(\hat G)=[(\omega\otimes\iota)(W)\mid\omega\in K^*]\ \ \text{and}\ \ W\in M(C_0(G)\otimes C_0(\hat G)).
$$

The pentagon relation for $W$ can be written in the following equivalent forms:
$$
(\Delta\otimes\iota)(W)=W_{13}W_{23},\ \ (\iota\otimes\Dhat)(W)=W_{13}W_{12}.
$$

Denote by $J$, resp. $\hat J$, the modular involutions on $L^2(G)$ defined by $\varphi$, resp. $\hat\varphi$. Then $J$ and $\hat J$ commute up to a scalar factor. The unitary antipode on $M$, resp. $\hat M$, is given by $R(x)=\hat Jx^*\hat J$, resp. $\hat R(a)=Ja^*J$. We have $(R\otimes\hat R)(W)=W$, that is,
$$
(\hat J\otimes J)W^*(\hat J\otimes J)=W.
$$

In addition to $W$ it is convenient to use another multiplicative unitary $V$ corresponding to the GNS-representation defined by the right Haar weight $\psi$. It is defined by
$$
V=(\hat J\otimes\hat J)\hat W(\hat J\otimes\hat J)\in L^\infty(\hat G)'\bar\otimes L^\infty(G),
$$
and we have
$$
\Delta(x)=V(x\otimes1)V^*\ \ \text{for}\ \ x\in L^\infty(G).
$$

A locally compact quantum group $G$ is called regular~\cite{BS,BSV}, if
$$
(K\otimes1)W(1\otimes K)\subset K\otimes K,
$$
or equivalently, $(1\otimes K)W(K\otimes1)\subset K\otimes K$. This is equivalent to several other conditions. In particular, if $G$ is regular then $[(K\otimes1)W(1\otimes K)]= K\otimes K$ and $$[(C_0(G)\otimes1)W(1\otimes C_0(\hat G))]=C_0(G)\otimes C_0(\hat G).$$
Recall that genuine locally compact groups, and therefore also their duals, are always regular.

\subsection{Actions on operator algebras}
A left action of a locally compact quantum group~$G$ on a von Neumann algebra $N$ is a normal unital injective $*$-homomorphism $\alpha\colon N\to L^\infty(G)\bar\otimes N$ such that $(\iota\otimes\alpha)\alpha=(\Delta\otimes\iota)\alpha$. A continuous left action of $G$ on a C$^*$-algebra $A$ is a nondegenerate injective $*$-homomorphism $\alpha\colon A\to M(C_0(G)\otimes A)$ such that $(\iota\otimes\alpha)\alpha=(\Delta\otimes\iota)\alpha$ and the following cancelation property holds:
$$
C_0(G)\otimes A=[(C_0(G)\otimes1)\alpha(A)].
$$

The following proposition is a small variation of results of Baaj, Skandalis and Vaes, see Propositions 5.7 and 5.8 in \cite{BSV}. We include a complete proof for convenience.

\begin{proposition}\label{pBSV}
Assume $G$ is a regular locally compact quantum group, $N$ is a von Neumann algebra and $\alpha\colon N\to L^\infty(G)\bar\otimes N$ is an action of $G$ on $N$. For a subspace $X\subset N$ define
$
X_\alpha=[(\omega\otimes\iota)\alpha(X)\mid\omega\in K^*]\subset N.
$
Then for any C$^*$-subalgebra $A\subset N$, if $A_\alpha\subset A$, then $A_\alpha$ is a C$^*$-algebra and $\alpha|_{A_\alpha}$ is a continuous action of $G$ on $A_\alpha$.
\end{proposition}

\bp Since $AA_\alpha\subset A$ by assumption, we have
\begin{align*}
A_\alpha &\supset[(\omega\otimes\iota)\alpha(A(\nu\otimes\iota)\alpha(A))\mid\omega,\nu\in K^*]\\
&=[(\nu\otimes\omega\otimes\iota)(\alpha(A)_{23}(\Delta\otimes\iota)\alpha(A))\mid\omega,\nu\in K^*]\\
&=[(\nu\otimes\omega\otimes\iota)(\alpha(A)_{23}V_{12}\alpha(A)_{13}V^*_{12})\mid\omega,\nu\in K^*]\\
&=[(\nu\otimes\omega\otimes\iota)(\alpha(A)_{23}V_{12}\alpha(A)_{13})\mid\omega,\nu\in K^*]\\
&=[(\nu\otimes\omega\otimes\iota)(\alpha(A)_{23}\alpha(A)_{13})\mid\omega,\nu\in K^*],
\end{align*}
where in the last step we used that $[(K\otimes1)V(1\otimes K)]=K\otimes K$ by regularity. We thus see that $A_\alpha A_\alpha\subset A_\alpha$. It is also clear that $A_\alpha$ is invariant under the $*$-operation. Thus $A_\alpha$ is a C$^*$-algebra.

\smallskip

In order to show that $\alpha|_{A_\alpha}$ is an action, observe first that $(X_\alpha)_\alpha=X_\alpha$ for any subspace $X\subset N$. Therefore replacing $A$ by $A_\alpha$ we may assume that $A=A_\alpha$. We then have
\begin{align*}
[\alpha(A)(C_0(G)\otimes1)]&=[\alpha((\omega\otimes\iota)\alpha(A))(C_0(G)\otimes1)\mid\omega\in K^*]\\
&=[(\omega\otimes\iota\otimes\iota)(V_{12}\alpha(A)_{13}V^*_{12}(1\otimes C_0(G)\otimes1))\mid\omega\in K^*]\\
&=[(\omega\otimes\iota\otimes\iota)(V_{12}\alpha(A)_{13}(1\otimes C_0(G)\otimes1))\mid\omega\in K^*]\\
&=[(\omega\otimes\iota\otimes\iota)(\alpha(A)_{13}(1\otimes C_0(G)\otimes1))\mid\omega\in K^*],
\end{align*}
where in the last step we used that $[(K\otimes1)V(1\otimes C_0(G)]=K\otimes C_0(G)$ by regularity. Therefore
$$
[\alpha(A)(C_0(G)\otimes1)]=C_0(G)\otimes A.
$$
From this we conclude that $\alpha(A)\subset M(C_0(G)\otimes A)$ and $\alpha|_A$ is a continuous action of $G$.
\ep

\subsection{Crossed products and duality} \label{ssCP}
Given a continuous left action $\alpha$ of a locally compact quantum group $G$ on a C$^*$-algebra $A$, the reduced C$^*$-crossed product $G\ltimes_\alpha A$ (since we are going to consider only reduced crossed products, we omit $r$ in the notation) is defined by
$$
G\ltimes_\alpha A=[(C_0(\hat G)\otimes1)\alpha(A)]\subset M(K\otimes A).
$$
It is equipped with the dual continuous right action of $\hat G$, or in other words, with a continuous left action $\hat\alpha$ of $\Ghop$, which is the quantum group $\hat G$ with the opposite comultiplication $\Dhop$ on $L^\infty(\Ghop)=L^\infty(\hat G)$. Namely, we have
$$
\hat\alpha(x)=(\Whop\otimes1)^*(1\otimes x)(\Whop\otimes1) \ \ \text{for}\ \ x\in G\ltimes_\alpha A,
$$
where
$$
\Whop=(J\otimes J)\hat W(J\otimes J)\in L^\infty(\hat G)\bar\otimes L^\infty(G)'
$$
is the multiplicative unitary of $\Ghop$ (see \cite[Section~4]{KVvN}), so
$$
\hat\alpha(\alpha(a))=1\otimes\alpha(a)
\ \ \text{for}\ \ a\in A\ \ \text{and}\ \ \hat\alpha(x\otimes1)=\Dhop(x)\otimes1\ \ \text{for}\ \ x\in C_0(\hat G).
$$
Then the double crossed product is
$$
\Ghop\ltimes_{\hat\alpha}G\ltimes_\alpha A
=[(JC_0(G)J\otimes1\otimes1)(\Dhop(C_0(\hat G))\otimes1)(1\otimes\alpha(A))].
$$
Since
$$
\Dhop(x)=W(x\otimes1)W^*\ \ \text{for}\ \ x\in C_0(\hat G),
$$
the map $\Ad(W^*\otimes1)$ maps $\Ghop\ltimes_{\hat\alpha}G\ltimes_\alpha A$ onto
$$
[(JC_0(G)JC_0(\hat G)\otimes1\otimes1)(\iota\otimes\alpha)\alpha(A)].
$$
In particular, if $[JC_0(G)JC_0(\hat G)]=K$, which is another equivalent formulation of regularity of $G$, we get the Takesaki-Takai duality
$$
\Ghop\ltimes_{\hat\alpha}G\ltimes_\alpha A\cong K\otimes\alpha(A)\cong K\otimes A.
$$

Assume now that $G$ is regular and $\beta$ is a continuous left action of $\Ghop$ on a C$^*$-algebra $B$. Assume also that there exists a unitary $X$ in $M(C_0(G)\otimes B)$ such that
$$
(\Delta\otimes\iota)(X)=X_{13}X_{23}\ \ \text{and}\ \ (\iota\otimes\beta)(X)=W_{12}X_{13}.
$$
Consider the $*$-homomorphism
$$
\eta\colon B\to M(K\otimes B),\ \ \eta(x)=X^*\beta(x)X.
$$
Then by a Landstad-type result of Vaes~\cite[Theorem~6.7]{V}, the space
$$
A=[(\omega\otimes\iota)\eta(B)\mid\omega\in K^*]\subset M(B)
$$
is a C$^*$-algebra, the formula
$$
\alpha(a)=X^*(1\otimes a)X
$$
defines a continuous left action of $G$ on $A$, and $\eta$ defines an isomorphism $B\cong G\ltimes_\alpha A$ intertwining~$\beta$ with~$\hat\alpha$. Note that if we already have $(B,\beta)=(G\ltimes_\alpha A,\hat\alpha)$, then we can take $X=W\otimes1$, in which case $\eta$ becomes the map defining the Takesaki-Takai isomorphism.

\bigskip

\section{Dual cocycles} \label{s2}

\subsection{Twisted group algebras}\label{ssTwisted}
Assume $G$ is a locally compact quantum group. By a measurable unitary dual $2$-cocycle on $G$, or a measurable unitary $2$-cocycle on $\hat G$, we mean a unitary element $\Omega\in L^\infty(\hat G)\bar\otimes L^\infty(\hat G)$ such that
$$
(\Omega\otimes1)(\Dhat\otimes\iota)(\Omega)=(1\otimes\Omega)(\iota\otimes\Dhat)(\Omega).
$$
We say that $\Omega$ is continuous if $\Omega\in M(C_0(\hat G)\otimes C_0(\hat G))$.

Given a measurable unitary $2$-cocycle $\Omega$, the cocycle condition can be written as
\begin{equation} \label{ecocycle0}
(\Dhat\otimes\iota)(\hat W\Omega^*)\Omega^*_{12}=(\hat W\Omega^*)_{13}(\hat W\Omega^*)_{23}.
\end{equation}
Indeed, we have
\begin{align*}
(\Dhat\otimes\iota)(\hat W\Omega^*)&=\hat W_{13}\hat W_{23}(\Dhat\otimes\iota)(\Omega^*)\\
&=\hat W_{13}\hat W_{23}(\iota\otimes\Dhat)(\Omega^*)(1\otimes\Omega^*)(\Omega\otimes1)\\
&=\hat W_{13}\hat W_{23}\hat W^*_{23}\Omega^*_{13}\hat W_{23}\Omega^*_{23}\Omega_{12},
\end{align*}
which is what we claimed.

Identity \eqref{ecocycle0} shows that the space of operators $(\omega\otimes\iota)(\hat W\Omega^*)$, $\omega\in K^*$, forms an algebra. The C$^*$-algebra $C^*_r(\hat G;\Omega)$ generated by this algebra is called the reduced twisted group C$^*$-algebra of~$\hat G$. The von Neumann algebra $C^*_r(\hat G;\Omega)''\subset B(L^2(G))$ is denoted by $W^*(\hat G;\Omega)$.

The following theorem in full generality is quite nontrivial and follows from results of De Commer~\cite[Propositions~11.2.1 and~11.2.2]{DCth}, which, in turn, rely on an analogue of manageability of multiplicative unitaries for measured quantum groupoids established by Enock~\cite{En}.

\begin{theorem} \label{tCdef}
We have $C^*_r(\hat G;\Omega)=[(\omega\otimes\iota)(\hat W\Omega^*)\mid \omega\in K^*]$ and $W\Omega^*\in M(K\otimes C^*_r(\hat G;\Omega))$.
\end{theorem}

For regular quantum groups the theorem is, however, not difficult to prove. Indeed, when $G$ is a compact quantum group, the equality $C^*_r(\hat G;\Omega)=[(\omega\otimes\iota)(\hat W\Omega^*)\mid \omega\in K^*]$ was proved in \cite[Lemma~4.9]{BDRV}. The same arguments work for any regular locally compact quantum group. In Section~\ref{ssDeform} we will also give a proof of this equality for arbitrary locally compact quantum groups that is independent of results of De Commer, by constructing a different set of generators of~$C^*_r(\hat G;\Omega)$. 

On the other hand, to show that $W\Omega^*\in M(K\otimes C^*_r(\hat G;\Omega))$ for regular quantum groups we can adapt the proof of \cite[Proposition~3.6]{BS} of a similar result for the multiplicative unitary. For this, rewrite identity~\eqref{ecocycle0} as
\begin{equation}\label{ecocycle2}
\hat W^*_{12}(\hat W\Omega^*)_{23}(\hat W\Omega^*)_{12}=(\hat W\Omega^*)_{13}(\hat W\Omega^*)_{23}.
\end{equation}
Multiplying by $K\otimes1\otimes1$ on the right and applying the slice maps to the second leg we get
$$
[(\iota\otimes\omega\otimes\iota)(\hat W^*_{12}(\hat W\Omega^*)_{23}(K\otimes1\otimes1))\mid\omega\in K^*]=\hat W\Omega^*(K\otimes C^*_r(\hat G;\Omega)).
$$
Using that  $[(1\otimes K)\hat W^*(K\otimes1)]=K\otimes K$ by regularity, we see that the left hand side equals $K\otimes C^*_r(\hat G;\Omega)$, so
$$
\hat W\Omega^*(K\otimes C^*_r(\hat G;\Omega))=K\otimes C^*_r(\hat G;\Omega).
$$
Similarly, rewriting \eqref{ecocycle2} as
$$
(\hat W\Omega^*)_{12}(\hat W\Omega^*)^*_{23}=(\hat W\Omega^*)^*_{23}\hat W_{12}(\hat W\Omega^*)_{13},
$$
multiplying this identity by $K\otimes1\otimes1$ on the left and applying the slice maps to the second leg, we get
$$
K\otimes C^*_r(\hat G;\Omega)=(K\otimes C^*_r(\hat G;\Omega))\hat W\Omega^*.
$$
Therefore $\hat W\Omega^*\in M(K\otimes C^*_r(\hat G;\Omega))$.

\smallskip

Let us also note the following.

\begin{proposition}
If $\Omega$ is a continuous unitary $2$-cocycle on $\hat G$, then $$\hat W\Omega^*\in M(C_0(\hat G)\otimes C^*_r(G;\Omega)).$$
\end{proposition}

\bp Using identity \eqref{ecocycle2} in the form
$$
(\hat W\Omega^*)_{13}=\hat W^*_{12}(\hat W\Omega^*)_{23}\hat W_{12}\Omega^*_{12}(\hat W\Omega^*)_{23}^*,
$$
from $\hat W\Omega^*\in M(K\otimes C^*_r(\hat G;\Omega))$ we see that $(\hat W\Omega^*)_{13}\in M(C_0(\hat G)\otimes K\otimes C^*_r(\hat G;\Omega))$. Applying the slice maps to the second leg we conclude that
$\hat W\Omega^*\in M(C_0(\hat G)\otimes C^*_r(G;\Omega))$.
\ep

The von Neumann algebras $W^*(\hat G;\Omega)$ (in fact, more general von Neumann-algebraic cocycle crossed products) were extensively studied by Vaes and Vainerman~\cite{VV}. In particular, in \cite[Proposition~1.4]{VV} they showed that there exists a right action $\beta$ of $G$ on $W^*(\hat G;\Omega)$ such that
$$
(\iota\otimes\beta)(\hat W\Omega^*)=\hat W_{13}(\hat W\Omega^*)_{12}.
$$
This action is given by
$$
\beta(x)=V(x\otimes1)V^*\ \ \text{for}\ \ x\in W^*(\hat G;\Omega).
$$
Another useful formula, which follows from \eqref{ecocycle2}, see \cite[Proposition 1.5]{VV}, is
\begin{equation} \label{eaction0}
\beta(x)=(\hat W\Omega^*)_{21}(1\otimes x)(\hat W\Omega^*)_{21}^*.
\end{equation}

\begin{proposition} \label{paction1}
The restriction of $\beta$ to $C^*_r(\hat G;\Omega)$ defines a continuous action of $G$ on the C$^*$-algebra $C^*_r(\hat G;\Omega)$.
\end{proposition}

\bp Since $\hat W\in M(K\otimes C_0(G))$, from the equality $(\iota\otimes\beta)(\hat W\Omega^*)=\hat W_{13}(\hat W\Omega^*)_{12}$ we get
$$
(K\otimes1\otimes C_0(G))(\iota\otimes\beta)(\hat W\Omega^*) =(K\otimes1\otimes C_0(G))(\hat W\Omega^*)_{12}.
$$
Applying the slice maps to the first leg we get
$$
[(1\otimes C_0(G))\beta(C^*_r(\hat G;\Omega))]=C^*_r(\hat G;\Omega)\otimes C_0(G),
$$
which proves the proposition.
\ep

\subsection{Deformed quantum group} \label{ssDG}
Given a unitary 2-cocycle $\Omega\in L^\infty(\hat G)\bar\otimes L^\infty(\hat G)$, we can define a new coproduct $\Dhat_\Omega$ on $L^\infty(\hat G)$ by
$$
\Dhat_\Omega(x)=\Omega\Dhat(x)\Omega^*\ \ \text{for}\ \ x\in L^\infty(\hat G).
$$
By a result of De Commer \cite{DC}, the pair $\hat G_\Omega=(L^\infty(\hat G),\Dhat_\Omega)$ is again a locally compact quantum group. We will use the subscript $\Omega$ to denote the objects related to $\hat G_\Omega$, such as the coproduct, the multiplicative unitary, etc.

In order to describe the multiplicative unitary $\hat W_\Omega$ of $\hat G_\Omega$ we need to recall some results of Vaes and Vainerman \cite{VV}. By \cite[Lemma 1.12]{VV} the action $\beta$ of $G$ on $W^*(\hat G;\Omega)$ is integrable, meaning that $(\iota\otimes\varphi)\beta$ is a n.s.f. operator valued weight from $W^*(\hat G;\Omega)$ to
$W^*(\hat G;\Omega)^\beta=\C1$. Therefore we have a n.s.f. weight $\tilde\varphi$ on $W^*(\hat G;\Omega)$ such that
$$
 \tilde\varphi(x)1=(\iota\otimes\varphi)\beta(x)\ \ \text{for} \ \ x\in W^*(\hat G;\Omega)_+.
$$
By construction $W^*(\hat G;\Omega)$ is represented on $L^2(G)$. By \cite[Proposition 1.15]{VV} this representation can be identified with the GNS-representation defined by the weight $\tilde\varphi$, with the corresponding map $\tilde\Lambda\colon\mathcal N_{\tilde\varphi}\to L^2(G)$ given by
\begin{equation} \label{etildelambda}
\tilde\Lambda((\omega\otimes\iota)(\hat W\Omega^*))=\Lambda((\omega\otimes\iota)(\hat W))
\end{equation}
for suitable $\omega\in K^*$. Denote by $\tilde J$ the modular involution on $L^2(G)$ defined by $\tilde\varphi$. The von Neumann algebra $L^\infty(\hat G_\Omega)=L^\infty(\hat G)\subset B(L^2(G))$ is in the standard form, so $\hat J_\Omega=\hat J$, and by \cite[Proposition~5.4]{DC} we have
\begin{equation}\label{eDC1}
\hat W_\Omega=(\tilde J\otimes \hat J)\Omega\hat W^*(J\otimes\hat J)\Omega^*.
\end{equation}
From this we immediately get the following proposition.

\begin{proposition} \label{pcod}
For any measurable unitary $2$-cocycle $\Omega\in L^\infty(\hat G)\bar\otimes L^\infty(\hat G)$ on $\hat G$, the element $\Omega^*$ is a measurable unitary $2$-cocycle on $\hat G_\Omega$, and we have
$$
C^*_r(\hat G_\Omega;\Omega^*)=\hat J C^*_r(\hat G;\Omega)\hat J.
$$
\end{proposition}

Therefore $C^*_r(\hat G;\Omega)$ is $*$-anti-isomorphic to $C^*_r(\hat G_\Omega;\Omega^*)$. By Proposition~~\ref{paction1} we have a continuous right action of $G_\Omega$ on $C^*_r(\hat G_\Omega;\Omega^*)$, where $G_\Omega$ is the dual of $\hat G_\Omega$. Using the unitary antipode $R_\Omega(x)=\hat J_\Omega x^*\hat J_\Omega=\hat J x^*\hat J$ on~$C_0(G_\Omega)$ we can transform this action to a continuous left action $\beta_\Omega$ of~$G_\Omega$ on~$C^*_r(\hat G;\Omega)$.

\begin{lemma}\label{laction}
We have
$$
\beta_\Omega(x)=W^*_\Omega(1\otimes x)W_\Omega\ \ \text{for}\ \ x\in C^*_r(\hat G;\Omega),
$$
and
$$
(\iota\otimes\beta_\Omega)(\hat W\Omega^*)=(\hat W\Omega^*)_{13}\big((\tilde J\otimes\hat J)\hat W^*_\Omega(\tilde J\otimes\hat J)\big)_{12}.
$$
\end{lemma}

\bp The right action of $G_\Omega$ on $C^*_r(\hat G_\Omega;\Omega^*)$ is given by $x\mapsto V_\Omega(x\otimes1)V_\Omega^*$. Therefore the left action of $G_\Omega$ on $C^*_r(\hat G;\Omega)$ is defined by
$$
\beta_\Omega(x)=(\hat J\otimes\hat J)(V_\Omega(\hat Jx\hat J\otimes1)V_\Omega^*)_{21}(\hat J\otimes\hat J).
$$
Since
$$
(V_\Omega)_{21}=(\hat J\otimes\hat J)(\hat W_\Omega)_{21}(\hat J\otimes\hat J)=(\hat J\otimes\hat J)W_\Omega^*(\hat J\otimes\hat J),
$$
we get the first formula for $\beta_\Omega$ in the formulation.

Similarly, since the right action of $G_\Omega$ on $C^*_r(\hat G_\Omega;\Omega^*)$ maps $(\omega\otimes\iota)(\hat W_\Omega\Omega)$ into $$(\omega\otimes\iota\otimes\iota)((\hat W_\Omega)_{13}(\hat W_\Omega\Omega)_{12}),$$
we have
$$
(\omega\otimes\beta_\Omega)((J_1\otimes\hat J)\hat W_\Omega\Omega(J_2\otimes\hat J))=(\omega\otimes\iota\otimes\iota)((J_1\otimes\hat J\otimes\hat J)
(\hat W_\Omega)_{12}(\hat W_\Omega\Omega)_{13}(J_2\otimes\hat J\otimes\hat J))
$$
for any $\omega\in K^*$ and any bounded antilinear operators $J_1$ and $J_2$. Since
$$
\hat W_\Omega\Omega=(\tilde J\otimes \hat J)(\hat W\Omega^*)^*(J\otimes\hat J),
$$
taking $J_1=\tilde J$ and $J_2=J$ we get
$$
(\iota\otimes\beta_\Omega)((\hat W\Omega^*)^*)=\big((\tilde J\otimes\hat J)\hat W_\Omega(\tilde J\otimes\hat J)\big)_{12}(\hat W\Omega^*)^*_{13},
$$
which is exactly the second formula in the formulation.
\ep

Therefore we have a left action $\beta_\Omega$ of $G_\Omega$ and a right action $\beta$ of $G$ on~$C^*_r(\hat G;\Omega)$. Using that $(\iota\otimes\beta)(\hat W\Omega^*)=\hat W_{13}(\hat W\Omega^*)_{12}$ and the second formula in the lemma above, we see that these actions commute: $(\beta_\Omega\otimes\iota)\beta=(\iota\otimes\beta)\beta_\Omega$.

\subsection{Twisted crossed products}\label{ssTCP}
Assume $\Omega$ is a measurable unitary $2$-cocycle on $\hat G$ and $\alpha$ is a continuous left action of $\Ghop$ on a C$^*$-algebra $A$.

\begin{definition}
The {\em reduced twisted crossed product} $\Ghop\ltimes_{\alpha,\Omega}A$ is defined as the C$^*$-subalgebra of~$M(K\otimes A)$ generated by
$
(J\hat JC^*_r(\hat G;\Omega)\hat JJ\otimes1)\alpha(A).
$
\end{definition}

\begin{proposition}\label{pdeftw}
We have $\Ghop\ltimes_{\alpha,\Omega}A=[(J\hat JC^*_r(\hat G;\Omega)\hat JJ\otimes1)\alpha(A)]$, and the formula
$$
\hat\alpha(x)=\Ad\big((1\otimes J\hat J\otimes1)(W^*_\Omega\otimes1)(1\otimes\hat J J\otimes1)\big)(1\otimes x)
$$
defines a continuous left action of $G_\Omega$ on $\Ghop\ltimes_{\alpha,\Omega}A$.
\end{proposition}

\bp The first part is proved in the standard way. Namely, observe first that
$$
\Whop=(J\otimes J)\hat W(J\otimes J)=(1\otimes J\hat J)\hat W^*(1\otimes\hat JJ),
$$
whence
$$
\Dhop(x)=\Ad\big((1\otimes J\hat J)\hat W\Omega^*(1\otimes\hat JJ)\big)(1\otimes x)\ \ \text{for}\ \ x\in L^\infty(\hat G).
$$
It follows that
\begin{multline*}
(1\otimes J\hat J\otimes1)(\hat W\Omega^*\otimes1)(1\otimes\hat JJ\otimes1)(1\otimes\alpha(A))\\
=(\iota\otimes\alpha)\alpha(A)(1\otimes J\hat J\otimes1)(\hat W\Omega^*\otimes1)(1\otimes\hat JJ\otimes1).
\end{multline*}
Applying the slice maps to the first leg and using that $(K\otimes1)\alpha(A)\subset K\otimes A$ we conclude that
$$
[(J\hat JC^*_r(\hat G;\Omega)\hat JJ\otimes1)\alpha(A)]\subset
[\alpha(A)(J\hat JC^*_r(\hat G;\Omega)\hat JJ\otimes1)],
$$
and therefore
$$
\Ghop\ltimes_{\alpha,\Omega}A=[\alpha(A)(J\hat JC^*_r(\hat G;\Omega)\hat JJ\otimes1)]
=[(J\hat JC^*_r(\hat G;\Omega)\hat JJ\otimes1)\alpha(A)].
$$

For the second part, note that the element
$$
(1\otimes J\hat J\otimes1)(W^*_\Omega\otimes1)(1\otimes\hat J J\otimes1)
\in L^\infty(G_\Omega)\bar\otimes L^\infty(\hat G)'\otimes1
$$
commutes with $1\otimes\alpha(A)$. Therefore it suffices to check that the formula
$$
J\hat JC^*_r(\hat G;\Omega)\hat JJ\ni x\mapsto \Ad\big((1\otimes J\hat J)W^*_\Omega(1\otimes\hat J J)\big)(1\otimes x)
$$
defines a continuous action of $G_\Omega$ on $J\hat JC^*_r(\hat G;\Omega)\hat JJ$. But this is true by Lemma~\ref{laction}.
\ep

Twisted crossed products, or cocycle crossed products, in the von Neumann algebraic setting were defined by Vaes and Vainerman~\cite{VV}. Our definition is of course related to theirs, but not in the most straightforward way. Namely, assume we are given a left action $\alpha$ of $\Ghop$ on a von Neumann algebra $N$. Then $(\alpha,\Omega_{21}^*)$ is a cocycle action of $(\Ghop)_{\Omega_{21}}=\Ghopo$ on $N$ in the sense of~\cite[Definition~1.1]{VV}. The von Neumann-algebraic cocycle crossed product of $N$ by $(\Ghop)_{\Omega_{21}}$ is defined as the von Neumann algebra generated by $\alpha(N)$ and $W^*((\Ghop)_{\Omega_{21}};\Omega_{21}^*)\otimes1$, see~\cite[Definition~1.3]{VV}.

\begin{lemma} \label{ltwVV}
Letting $X=\tilde JJ$, we have $X\in L^\infty(\hat G)$ and
$$
J\hat JW^*(\hat G;\Omega)\hat JJ=\hat JX^* \hat JW^*((\Ghop)_{\Omega_{21}};\Omega_{21}^*)\hat JX\hat J.
$$
\end{lemma}

\bp The claim that $X\in L^\infty(\hat G)$ is in~\cite[Section~5]{DC}. By Proposition~\ref{pcod}, applied to the cocycle~$\Omega_{21}$ on $\Ghop$, we have
$$
W^*((\Ghop)_{\Omega_{21}};\Omega_{21}^*)=\hat JW^*(\Ghop;\Omega_{21})\hat J.
$$
By \cite[Proposition~6.3]{DC} we also have
$$
\Omega^*(X\otimes X)=\Dhat(X)(\hat R\otimes\hat R)(\Omega_{21})=\Dhat(X)(J\otimes J)\Omega_{21}^*(J\otimes J),
$$
whence
\begin{equation} \label{eid}
\hat W\Omega^*(X\otimes X)=(1\otimes X)\hat W(J\otimes J)\Omega_{21}^*(J\otimes J)
=(1\otimes X)(J\otimes J)\Whop\Omega_{21}^*(J\otimes J),
\end{equation}
and therefore
$$
W^*(\hat G;\Omega)=XJW^*(\Ghop;\Omega_{21})JX^*=\tilde JW^*(\Ghop;\Omega_{21})\tilde J.
$$
It follows that
$$
J\hat JW^*(\hat G;\Omega)\hat JJ=\hat J J\tilde JW^*(\Ghop;\Omega_{21})\tilde J J\hat J=\hat J J\tilde J\hat JW^*((\Ghop)_{\Omega_{21}};\Omega_{21}^*)\hat J\tilde JJ\hat J,
$$
which is what we need.
\ep

Therefore up to conjugation by $\hat JX^*\hat J\otimes1$ our definition of the twisted crossed product by $\Ghop$ is a C$^*$-algebraic version of the definition of Vaes and Vainerman of the cocycle crossed product by~$(\Ghop)_{\Omega_{21}}$.

\subsection{Regular cocycles}\label{ssRegular} Assume $\Omega$ is a measurable unitary cocycle on $\hat G$. Recall that we have a continuous right action $\beta$ of $G$ on $C^*_r(\hat G;\Omega)$, so we can consider the reduced crossed product
$$
C^*_r(\hat G;\Omega)\rtimes_\beta G=[\beta(C^*_r(\hat G;\Omega))(1\otimes\hat JC_0(\hat G)\hat J)].
$$
By \eqref{eaction0} we have $\beta(x)=(\hat W\Omega^*)_{21}(1\otimes x)(\hat W\Omega^*)_{21}^*$. The unitary $(\hat W\Omega^*)^*_{21}$ commutes with $1\otimes\hat JC_0(\hat G)\hat J$. Therefore the conjugation by this unitary maps $C^*_r(\hat G;\Omega)\rtimes_\beta G$ onto $1\otimes[C^*_r(\hat G;\Omega)\hat JC_0(\hat G)\hat J]$.

\begin{definition}
A cocycle $\Omega$ is called {\em regular} if $[C^*_r(\hat G;\Omega)\hat JC_0(\hat G)\hat J]=K.$
\end{definition}

Note that by a version of the Takesaki duality \cite[Proposition~1.20]{VV} the von Neumann algebra generated by $C^*_r(\hat G;\Omega)\hat JC_0(\hat G)\hat J$ coincides with $B(L^2(G))$ (this will also become clear from the proof of Proposition~\ref{p2cond} below). Therefore regularity of $\Omega$ is equivalent to the formally weaker condition
$$
C^*_r(\hat G;\Omega)\hat JC_0(\hat G)\hat J\subset K.
$$
Since the representation of $C^*_r(\hat G;\Omega)\rtimes_\beta G$ on $L^2(G)$ is faithful and irreducible, yet another equivalent formulation of regularity of $\Omega$ is that the C$^*$-algebra $C^*_r(\hat G;\Omega)\rtimes_\beta G$ is isomorphic to the algebra of compact operators on some Hilbert space.

By definition, regularity of the trivial cocycle $1$ is the same as regularity of $G$. Therefore regularity of cocycles is definitely not automatic. Even for regular locally compact quantum groups regularity of a cocycle is a very delicate question.  The only easy cases seem to be covered by the following proposition.

\begin{proposition}
Any measurable unitary $2$-cocycle on $\hat G$ is regular in the following cases:
\enu{i} $\hat G$ is a genuine locally compact group;
\enu{ii} $\hat G$ is a discrete quantum group.
\end{proposition}

\bp Part (i) is well-known and is proved in the same way as regularity of $\hat G$, by observing that the space $C^*_r(\hat G;\Omega)\hat JC_0(\hat G)\hat J$ contains a lot of integral operators. Part (ii) is obvious, as already the algebra $C_0(\hat G)$ consists of compact operators.
\ep

In view of various equivalent characterizations of regularity of quantum groups, it is natural to wonder how regularity of a cocycle is related to properties like $(K\otimes1)\hat W\Omega^*(1\otimes K)\subset K\otimes K$. We have the following result.

\begin{proposition} \label{p2cond}
For a cocycle $\Omega$ on $\hat G$ consider the following conditions:
\enu{i} $\Omega$ is regular;
\enu{ii} $(K\otimes1)\hat W\Omega^*(1\otimes K)\subset K\otimes K$.

\smallskip

Then $(\text{\rm i})\Rightarrow(\text{\rm ii})$. If $G$ is regular, then the two conditions are equivalent.
\end{proposition}

\bp In this proof it will be convenient to consider the right action $\beta$ of $G$ on $C^*_r(\hat G;\Omega)$ as the left action $\beta'$ of $\Gop$, so $\beta'(x)=V_{21}(1\otimes x)V^*_{21}$ for $x\in C^*_r(\hat G;\Omega)$. On the von Neumann algebra level, up to stabilization this action is dual. Namely, by \cite[Propositions~1.8 and~1.9]{VV} the unitary $Y=\hat V^*_{21}\Omega^*_{21}$ defines a left action $\gamma$ of $\hat G$ on $B(L^2(G))$ by
$$
\gamma(x)=Y(1\otimes x)Y^*,
$$
and we have an isomorphism
$$
W^*(\hat G;\Omega)\bar\otimes B(L^2(G))\cong (L^\infty(G)\otimes1\cup\gamma(B(L^2(G))))'',\ \ x\mapsto YxY^*,
$$
intertwining $\beta'\otimes\iota$ with the dual action $\hat\gamma$, defined by $\hat\gamma(x)=V_{21}(1\otimes x)V^*_{21}$.  Note also that
$$
\hat V^*_{21}=(J\otimes J)\hat W(J\otimes J)=\Whop
$$
and by~\eqref{eid},
$$
\hat W\Omega^*(\tilde J\otimes\tilde J)=(J\otimes\tilde J)\Whop\Omega^*_{21},
$$
so that
$$
Y=\hat V^*_{21}\Omega^*_{21}=(J\otimes \tilde J)\hat W\Omega^*(\tilde J\otimes\tilde J).
$$

We are now ready to prove the proposition.

\smallskip

Assume condition (i) holds. We claim that the restriction of $\gamma$ to $K$ defines a continuous action of~$\hat G$ on~$K$. Since $K=[C^*_r(\hat G;\Omega)\hat JC_0(\hat G)\hat J]$ and $Y$ commutes with $1\otimes C^*_r(\hat G;\Omega)$, it suffices to show that the restriction of $\gamma$ to $\hat JC_0(\hat G)\hat J$ defines a continuous action. But this is clear, since $\Omega_{21}$ commutes with $1\otimes\hat JC_0(\hat G)\hat J$ and therefore
$$
\gamma(\hat Jx\hat J)=\Whop(1\otimes\hat Jx\hat J)(\Whop)^*=
(J\otimes\hat J)(\Whop)^*(1\otimes x)\Whop(J\otimes\hat J)=
(J\otimes\hat J)\Dhop(x)(J\otimes\hat J).
$$

It follows that
$$
(K\otimes1)Y(1\otimes K)Y^*(K\otimes1)=(K\otimes1)\gamma(K)(K\otimes1)\subset K\otimes K.
$$
This means exactly that $(K\otimes1)Y(1\otimes K)\subset K\otimes K$, which is equivalent to (ii).

\smallskip

Assume now that $G$ is regular and condition (ii) holds. We claim again that the restriction of~$\gamma$ to~$K$ defines a continuous action of $\hat G$. By Proposition~\ref{pBSV} it suffices to show that $K_\gamma=K$. Condition~(ii) implies that $(K\otimes1)\gamma(K)(K\otimes1)\subset K\otimes K$, whence $K_\gamma\subset K$. By Proposition~\ref{pBSV} this already shows that $K_\gamma$ is a C$^*$-algebra. Since it is $\sigma$-strongly$^*$ dense in
$$
[(\omega\otimes\iota)\gamma(B(L^2(G)))\mid\omega\in K^*]''=B(L^2(G)),
$$
it follows that $K_\gamma=K$.

Next, we have an isomorphism
$$
C^*_r(\hat G;\Omega)\otimes K\cong\hat G\ltimes_\gamma K,\ \ x\mapsto YxY^*,
$$
intertwining $\beta'\otimes\iota$ with $\hat\gamma$. This is a C$^*$-algebraic version of \cite[Propositions~1.8]{VV}, and the proof is basically the same. Briefly, we have the identity
$$
Y^*_{23}\hat W_{12}Y_{23}=(\hat W\Omega^*)_{12}Y^*_{13},
$$
which is proved similarly to \eqref{ecocycle0}. Applying the slice maps to the first leg we conclude that $\Ad Y^*$ maps
$$
\hat G\ltimes_\gamma K=[(C_0(G)\otimes1)\gamma(K)]=[(\omega\otimes\iota\otimes\iota)(\hat W\otimes1)(1\otimes Y)(1\otimes1\otimes K)(1\otimes Y^*)\mid\omega\in K^*]
$$
onto
$$
[(\omega\otimes\iota\otimes\iota)(\hat W\Omega^*\otimes1)Y^*_{13}(1\otimes1\otimes K)\mid\omega\in K^*]=C^*_r(\hat G;\Omega)\otimes K,
$$
as claimed.

Consider now the double crossed product $\Gop\ltimes_{\hat\gamma}\hat G\ltimes_\gamma K$. By the Takesaki-Takai duality it is isomorphic to $K\otimes K$. What is however important to us, is only the equality
\begin{equation}\label{eid2}
[(\hat JC_0(\hat G)\hat J\otimes1)(\hat G\ltimes_\gamma K)]=K\otimes K,
\end{equation}
which is an immediate consequence of regularity of $G$, since $[\hat JC_0(\hat G)\hat J C_0(G)]=K$ and
$$
[(K\otimes1)\gamma(K)]=[(KC_0(\hat G)\otimes1)\gamma(K)]=K\otimes K.
$$
Applying~$\Ad Y^*$ to both sides of \eqref{eid2} and using that~$Y$ commutes with $\hat JC_0(\hat G)\hat J\otimes1$ we conclude that
$$
[\hat JC_0(\hat G)\hat JC^*_r(\hat G;\Omega)]\otimes K=K\otimes K,
$$
so $\Omega$ is regular.
\ep

\bigskip

\section{Deformation of C$^*$-algebras}\label{s3}

\subsection{Quantization maps}\label{squantmaps} Let $\Omega$ be a measurable unitary $2$-cocycle on a locally compact quantum group~$\hat G$.

\begin{proposition}\label{pquantmaps}
For every $\nu\in K^*$, the formula
$$
T_\nu(x)=(\iota\otimes\nu)(\hat W_\Omega\Omega(x\otimes1)(\hat W_\Omega\Omega)^*)
$$
defines a right $G$-equivariant map $ C_0(G)\to C^*_r(\hat G;\Omega)$, where we consider the right actions $\Delta$ and $\beta$ of~$G$ on $C_0(G)$ and $C^*_r(\hat G;\Omega)$, respectively. Furthermore, we have
$$
[T_\nu(C_0(G))\mid\nu\in K^*]=[(\omega\otimes\iota)(\hat W\Omega^*)\mid \omega\in K^*]=C^*_r(\hat G;\Omega).
$$
\end{proposition}

\bp
The proof relies on the following identity:
\begin{equation} \label{ecocycle4}
(\hat W_\Omega\Omega)_{23}\hat W_{12}(\hat W_\Omega\Omega)^*_{23}=(\hat W\Omega^*)_{12}(\hat W_\Omega\Omega)_{13}.
\end{equation}
To prove it, write identity \eqref{ecocycle2} for the cocycle $\Omega^*$ on $\hat G_\Omega$ as
$$
(\hat W_\Omega\Omega)_{23}(\hat W_\Omega\Omega)_{12}(\hat W_\Omega\Omega)^*_{23}=(\hat W_\Omega)_{12}(\hat W_\Omega\Omega)_{13}.
$$
Substituting $(\hat W_\Omega)_{12}$ on both sides of the above identity with $\big((\tilde J\otimes \hat J)\Omega(J\otimes \hat J)\hat W\Omega^*\big)_{12}$ we get
$$
(\hat W_\Omega\Omega)_{23}\big((\tilde J\otimes \hat J)\Omega(J\otimes \hat J)\big)_{12}\hat W_{12}(\hat W_\Omega\Omega)^*_{23}=\big((\tilde J\otimes \hat J)\Omega(J\otimes \hat J)\big)_{12}(\hat W\Omega^*)_{12}(\hat W_\Omega\Omega)_{13}.
$$
Since $(\hat W_\Omega\Omega)_{23}$ and $((\tilde J\otimes \hat J)\Omega(J\otimes \hat J))_{12}$ commute, this is exactly \eqref{ecocycle4}.

Applying the slice maps to the first and the third legs of~\eqref{ecocycle4} we see that the image of $T_\nu$ is contained in $C^*_r(\hat G;\Omega)$ and
$$
[T_\nu(C_0(G))\mid\nu\in K^*]=[(\omega\otimes\iota)(\hat W\Omega^*)\mid \omega\in K^*].
$$

It remains to check $G$-equivariance. For $x\in C_0(G)$ we compute:
\begin{align*}
\beta(T_\nu(x))&=(\iota\otimes\iota\otimes\nu)(V_{12}(\hat W_\Omega\Omega)_{13}(x\otimes1\otimes1)(\hat W_\Omega\Omega)^*_{13}V^*_{12})\\
&=(\iota\otimes\iota\otimes\nu)((\hat W_\Omega\Omega)_{13}V_{12}(x\otimes1\otimes1)V_{12}^*(\hat W_\Omega\Omega)^*_{13})\\
&=(T_\nu\otimes\iota)\Delta(x),
\end{align*}
which finishes the proof of the proposition.
\ep

As a byproduct we get an alternative proof of part of Theorem~\ref{tCdef}, as promised earlier: the space $[T_\nu(C_0(G))\mid\nu\in K^*]$ is clearly self-adjoint, hence the algebra $[(\omega\otimes\iota)(\hat W\Omega^*)\mid \omega\in K^*]$ is a C$^*$-algebra, so it coincides with $C^*_r(\hat G;\Omega)$.

\smallskip

The map $T_\nu$ depends only on the restriction of $\nu$ to $W^*(\hat G_\Omega;\Omega^*)=\hat JW^*(\hat G;\Omega)\hat J$. It extends to a normal map $L^\infty(G)\to W^*(\hat G;\Omega)$, which we continue to denote by $T_\nu$. Note also that, since $\hat W\Omega^*\in M(K\otimes C^*_r(\hat G;\Omega))$ and $\hat W_\Omega\Omega\in M(K\otimes C^*_r(\hat G_\Omega;\Omega^*))$ by Theorem~\ref{tCdef}, identity \eqref{ecocycle4} implies that
$\hat W_\Omega\Omega(C_0(G)\otimes1)(\hat W_\Omega\Omega)^*$ is a nondegenerate C$^*$-subalgebra of
$$
M(C^*_r(\hat G;\Omega)\otimes C^*_r(\hat G_\Omega;\Omega^*))=M(C^*_r(\hat G;\Omega)\otimes\hat J C^*_r(\hat G;\Omega)\hat J).
$$
This implies that $T_\nu$ maps $M(C_0(G))$ into $M(C^*_r(\hat G;\Omega))$, and the map $T_\nu\colon M(C_0(G))\to M(C^*_r(\hat G;\Omega))$ is strictly continuous on bounded sets.

\begin{example}
Assume $G$ is the dual of a discrete group $\Gamma$, so $L^\infty(G)=W^*(\Gamma)\subset B(\ell^2(\Gamma))$ and $\Delta(\lambda_s)=\lambda_s\otimes\lambda_s$ for $s\in\Gamma$. Then $L^\infty(\hat G)=\ell^\infty(\Gamma)$, and a $2$-cocycle on $\hat G$ is a $2$-cocycle $\Omega\colon\Gamma\times\Gamma\to\T$ on $\Gamma$ in the usual sense. The multiplicative unitary $\hat W$ is defined by $\hat W(\delta_s\otimes\delta_t)=\delta_s\otimes\delta_{st}$. The twisted group C$^*$-algebra $C^*_r(\hat G;\Omega)$ is generated by the operators $\lambda^\Omega_s=\lambda_s\overline{\Omega(s,\cdot)}$ satisfying $\lambda^\Omega_{st}=\Omega(s,t)\lambda^\Omega_s\lambda^\Omega_t$. In this case we have $G_\Omega=G$, and \eqref{ecocycle4} gives us the known identity
$$
\hat W\Omega(\lambda_s\otimes1)(\hat W\Omega)^*=\lambda^\Omega_s\otimes\lambda^{\bar\Omega}_s.
$$
Therefore the maps $T_\nu\colon C^*_r(\Gamma)\to C^*_r(\Gamma;\Omega)$ are given by $T_\nu(\lambda_s)=\nu(\lambda^{\bar\Omega}_s)\lambda^\Omega_s$.
\ee
\end{example}

We call the maps $T_\nu$ the {\em quantization maps}. We will write $T^\Omega_\nu$ for $T_\nu$ when we want to stress that we consider the quantization maps defined by $\Omega$. We can also define {\em dequantization maps} going in the opposite direction, although we will not need them in this paper. Namely, for $\omega\in W^*(\hat G;\Omega)_*$ we can define
$$
S_\omega\colon C^*_r(\hat G;\Omega)\to C_0(G),\ \ S_\omega(x)=(\omega\otimes\iota)\beta(x).
$$
Since the elements of the form $T_\nu(x)$, with $x\in C_0(G)$ and $\nu\in W^*(\hat G_\Omega;\Omega^*)_*$, span a dense subspace of $C^*_r(\hat G;\Omega)$, the following computation shows that the image of $S_\omega$ is contained in $C_0(G)$:
$$
S_\omega T_\nu(x)=(\omega\otimes\iota)(T_\nu\otimes\iota)\Delta(x)
=(\nu\otimes\omega\otimes\iota)((W_\Omega\Omega)_{21}(1\otimes\Delta(x))(\hat W_\Omega\Omega)^*_{21})\in C_0(G).
$$
The maps $S_\omega$ are again right $G$-equivariant.

\subsection{$\Omega$-Deformation} \label{ssDeform}

Assume now that $A$ is a C$^*$-algebra and $\alpha$ is a continuous left action of $G$ on~$A$. Since, as we observed in the previous subsection, $\hat W_\Omega\Omega(C_0(G)\otimes1)(\hat W_\Omega\Omega)^*$ is a nondegenerate C$^*$-subalgebra of
$M(C^*_r(\hat G;\Omega)\otimes\hat J C^*_r(\hat G;\Omega)\hat J)$, the maps $T_\nu\otimes\iota\colon C_0(G)\otimes A\to C^*_r(\hat G;\Omega)\otimes A$ extend to maps
$$
T_\nu\otimes\iota \colon M(C_0(G)\otimes A)\to M(C^*_r(\hat G;\Omega)\otimes A),
$$
defined by $(T_\nu\otimes\iota)(x)=(\nu\otimes\iota\otimes\iota)((\hat W_\Omega\Omega)_{21}(1\otimes x)(\hat W_\Omega\Omega)_{21}^*)$.

\begin{definition}\label{domegadef}
The {\em $\Omega$-deformation} of a $A$ is the C$^*$-subalgebra
$$
A_\Omega\subset M(C^*_r(\hat G;\Omega)\otimes A)
$$
generated by elements of the form $(T_\nu\otimes\iota)\alpha(a)$ for all $\nu\in K^*$ and $a\in A$. The maps $(T_\nu\otimes\iota)\alpha\colon A\to A_\Omega$ are called the {\em quantization maps}.
\end{definition}

Note that since the maps $T_\nu$ are right $G$-equivariant, we immediately see that
$$
A_\Omega\subset\{x\in M(C^*_r(\hat G;\Omega)\otimes A)\mid (\beta\otimes\iota)(x)=(\iota\otimes\alpha)(x)\}.
$$

\smallskip

As a first example consider $A=C_0(G)$ with the action of $G$ on itself by left translations, so $\alpha=\Delta$. In this case, using that $(T_\nu\otimes\iota)\Delta(x)=\beta(T_\nu(x))$ for all $x\in C_0(G)$, we get
\begin{equation} \label{ese}
C_0(G)_\Omega=\beta(C^*_r(\hat G;\Omega))\cong C^*_r(\hat G;\Omega).
\end{equation}
This provides a different perspective on the action $\beta$ of $G$ on $W^*(\hat G;\Omega)$. This action was defined in~\cite{VV} as a dual action on a twisted crossed product. We can now say that $\beta$ is simply the right action of $G$ on itself that survives under deformation. More precisely, we have the following general result.

\begin{proposition}
Assume $A$ is a C$^*$-algebra equipped with a continuous left action $\alpha$ of a locally compact quantum group $G$ and a continuous right action $\gamma$ of a locally compact quantum group $H$ such that $(\iota\otimes\gamma)\alpha=(\alpha\otimes\iota)\gamma$. Then the restriction of $\iota\otimes\gamma\colon M(C^*_r(\hat G;\Omega)\otimes A)\to M(C^*_r(\hat G;\Omega)\otimes A\otimes C_0(H))$ to~$A_\Omega$ defines a continuous right action of $H$ on $A_\Omega$.
\end{proposition}

\bp For any $\nu\in K^*$ we have
\begin{align*}
[(1\otimes1\otimes C_0(H))(\iota\otimes\gamma)(T_\nu\otimes\iota)\alpha(A)]
&=[(T_\nu\otimes\iota\otimes\iota)(\alpha\otimes\iota)\big((1\otimes C_0(H))\gamma(A)\big)]\\
&=[(T_\nu\otimes\iota)\alpha(A)]\otimes C_0(H).
\end{align*}
This implies that
$$
[(1\otimes1\otimes C_0(H))(\iota\otimes\gamma)(A_\Omega)]=A_\Omega\otimes C_0(H).
$$
From this we conclude that $(\iota\otimes\gamma)(A_\Omega)\subset M(A_\Omega\otimes C_0(H))$ and the restriction of $\iota\otimes\gamma$ to $A_\Omega$ defines a continuous action of $H$.
\ep

\subsection{Relation to twisted crossed products}

Assume $B$ is a C$^*$-algebra and $\gamma$ is a continuous left action of $\Ghop$ on $B$. Consider the crossed product
$$
\Ghop\ltimes_\gamma B=[(JC_0(G)J\otimes1)\gamma(B)]
$$
and the dual action $\alpha=\hat\gamma$ of $G$ on $\Ghop\ltimes_\gamma B$ given by
$$
\alpha(x)=\Ad\big((1\otimes J\hat J\otimes1)(W^*\otimes1)(1\otimes \hat J J\otimes1)\big)(1\otimes x).
$$
Generalizing the isomorphism $C_0(G)_\Omega\cong C^*_r(\hat G;\Omega)$ we then have the following result.

\begin{proposition}\label{pdefvstw}
We have $(\Ghop\ltimes_\gamma B)_\Omega\cong \Ghop\ltimes_{\gamma,\Omega}B$. More precisely, $$\Ad\big((1\otimes J\hat J\otimes1)(\hat W\Omega^*)^*_{21}(1\otimes \hat J J\otimes1)\big)$$ maps the $\Omega$-deformation of $\Ghop\ltimes_\gamma B$ onto $1\otimes(\Ghop\ltimes_{\gamma,\Omega}B)$.
\end{proposition}

\bp The conjugation by $\hat JJ$ defines a left $G$-equivariant isomorphism $JC_0(G)J\cong C_0(G)$, where on $JC_0(G)J$ we consider the action given by
$$
x\mapsto \Ad\big((1\otimes J\hat J)W^*(1\otimes \hat J J)\big)(1\otimes x).
$$
Therefore by \eqref{ese} the $\Omega$-deformation of $JC_0(G)J$ is equal to
$$
(1\otimes J\hat J)\beta(C^*_r(\hat G;\Omega))(1\otimes \hat J J).
$$
From this, by definition of the $\Omega$-deformation, we conclude that $(\Ghop\ltimes_\gamma B)_\Omega$ is generated by
$$
(1\otimes J\hat J\otimes1)\beta(C^*_r(\hat G;\Omega))\otimes1)(1\otimes \hat J J\otimes1)(1\otimes\gamma(B)).
$$
Now recall that $\beta(x)=(\hat W\Omega^*)_{21}(1\otimes x)(\hat W\Omega^*)_{21}^*$ by \eqref{eaction0}. Observing also that the unitary $$(1\otimes J\hat J\otimes1)(\hat W\Omega^*)^*_{21}(1\otimes \hat J J\otimes1)$$ commutes with $1\otimes\gamma(B)$, we conclude that the conjugation by this unitary maps $(\Ghop\ltimes_\gamma B)_\Omega$ onto the C$^*$-algebra generated by
$
(1\otimes J\hat J C^*_r(\hat G;\Omega)\hat JJ\otimes1)(1\otimes\gamma(B)).
$
But this is exactly $1\otimes(\Ghop\ltimes_{\gamma,\Omega}B)$.
\ep

Turning to more general actions, recall that for regular quantum groups any action is stably exterior equivalent to a dual action. Therefore it is natural to expect that up to stabilization $\Omega$-deformations can be expressed in terms of twisted crossed products, at least under some regularity assumptions. In order to formulate the result recall that in Sections~\ref{ssTCP} and~\ref{ssRegular} we already used the unitaries
$$
X=\tilde JJ\in L^\infty(\hat G)\ \ \text{and}\ \ Y=\hat V^*_{21}\Omega^*_{21}=(J\otimes \tilde J)\hat W\Omega^*(\tilde J\otimes\tilde J).
$$
Using formula \eqref{eDC1} for $\hat W_\Omega$ we can also write
$$
Y=(1\otimes\tilde J\hat J)(\hat W_\Omega\Omega)^*(1\otimes\hat J\tilde J).
$$
We also define a map $\eta_\Omega\colon\alpha(A)\to M(\hat J C^*_r(\hat G;\Omega)\hat J\otimes C^*_r(\hat G;\Omega)\otimes A)$ by
\begin{equation}\label{eetaomega}
\eta_\Omega(\alpha(a))=(\hat W_\Omega\Omega)_{21}(1\otimes\alpha(a))(\hat W_\Omega\Omega)_{21}^*,
\end{equation}
so that $(T_\nu\otimes\iota)\alpha(a)=(\nu\otimes\iota\otimes\iota)\eta_\Omega(\alpha(a))$.

\begin{theorem}\label{ttwisted}
Assume $\Omega$ is a regular cocycle on a locally compact quantum group $\hat G$. Then for any C$^*$-algebra $A$ equipped with a continuous left action $\alpha$ of $G$ we have
$$
\Ghop\ltimes_{\hat\alpha,\Omega}G\ltimes_\alpha A\cong K\otimes A_\Omega.
$$
Explicitly, the map $$\Ad\big((\hat JJ\otimes1\otimes1)Y^*_{21}(J\hat J\otimes1\otimes1)\big)=\Ad\big((\hat JX^*\hat J\otimes1\otimes1)(\hat W_\Omega\Omega)_{21}(\hat J X\hat J\otimes1\otimes1)\big)$$ defines such an isomorphism. This map is trivial on $J\hat JC^*_r(\hat G;\Omega)\hat JJ\otimes 1\otimes 1$ and it maps $\Dhop(x)\otimes1\in\Dhop(C_0(\hat G))\otimes1$ into $x\otimes1\otimes1$ and $1\otimes\alpha(a)\in1\otimes\alpha(A)$ into $\Ad(\hat JX^*\hat J\otimes1\otimes1)\eta_\Omega(\alpha(a))$.
\end{theorem}

We will simultaneously prove the following.

\begin{theorem}\label{talgebra}
If $\Omega$ is regular, then
$$
A_\Omega=[(T_\nu\otimes\iota)\alpha(A)\mid\nu\in K^*].
$$
\end{theorem}

\bp[Proof of Theorems~\ref{ttwisted} and~\ref{talgebra}] From the equality $Y=(J\otimes \tilde J)\hat W\Omega^*(\tilde J\otimes\tilde J)$ it is clear that $Y$ commutes with $1\otimes C^*_r(\hat G;\Omega)$. Hence the homomorphism in the formulation is trivial on $J\hat JC^*_r(\hat G;\Omega)\hat JJ\otimes 1\otimes 1$. Clearly, it maps $1\otimes\alpha(a)\in1\otimes\alpha(A)$ into $\Ad(\hat JX^*\hat J\otimes1\otimes1)\eta_\Omega(\alpha(a))$. Finally, for $x\in C_0(\hat G)$ this homomorphism maps $\Dhop(x)\otimes1$ into
$$
(\hat W_\Omega\Omega)_{21}(\Dhop(x)\otimes1)(\hat W_\Omega\Omega^*)_{21}=(\hat W_\Omega\Omega\Dhat(x)(\hat W_\Omega\Omega)^*)_{21}
=(\hat W_\Omega\Dhat_\Omega(x)\hat W_\Omega^*)_{21}=x\otimes1\otimes1.
$$
Since $[J\hat JC^*_r(\hat G;\Omega)\hat JJC_0(\hat G)]=K$ by regularity of $\Omega$, it follows that this homomorphism maps the C$^*$-algebra $\Ghop\ltimes_{\hat\alpha,\Omega}G\ltimes_\alpha A$ onto
$$
[(K\otimes1\otimes1)\Ad(\hat JX^*\hat J\otimes1\otimes1)\eta_\Omega(\alpha(A))].
$$
In particular, the last space is a C$^*$-algebra, so it coincides with
$$
[(K\otimes1\otimes1)\Ad(\hat JX^*\hat J\otimes1\otimes1)\eta_\Omega(\alpha(A))(K\otimes1\otimes1)]=K\otimes[(\nu\otimes\iota\otimes\iota) \eta_\Omega(\alpha(A))\mid\nu\in K^*].
$$
This shows that $[(\nu\otimes\iota\otimes\iota) \eta_\Omega(\alpha(A))\mid\nu\in K^*]$ is a C$^*$-algebra and finishes the proof of both theorems.
\ep

Note that while the regularity of $\Omega$ is a necessary condition for the conclusion of Theorem~\ref{ttwisted} to be true, it is not clear whether this is the case for Theorem~\ref{talgebra}. For example, by the proof of Proposition~\ref{pdefvstw}, for dual actions Theorem~\ref{talgebra} remains true for any $\Omega$.

\subsection{Deformed action} Recall from Section~\ref{ssDG} that we have a continuous left action $\beta_\Omega$ of $G_\Omega$ on~$C^*_r(\hat G;\Omega)$, which commutes with the right action $\beta$ of $G$. This suggests that the action $\beta_\Omega\otimes\iota$ on $C^*_r(\hat G;\Omega)\otimes A$ defines a continuous action on $A_\Omega$. In other words, we want to define a left $\alpha_\Omega$ of $G_\Omega$ on $A_\Omega$ by
\begin{equation} \label{edefaction}
\alpha_\Omega(x)=(W^*_\Omega\otimes1)(1\otimes x)(W_\Omega\otimes1)\ \ \text{for}\ \ x\in A_\Omega\subset M(C^*_r(\hat G;\Omega)\otimes A).
\end{equation}
We can prove that this is indeed a continuous action of $G_\Omega$ under an additional regularity assumption.

\begin{theorem} \label{tdefaction}
Assume $\Omega$ is a measurable unitary $2$-cocycle on a locally compact quantum group $\hat G$ such that the deformed quantum group $G_\Omega$ is regular. Then for any C$^*$-algebra $A$ equipped with a continuous left action $\alpha$ of $G$, the formula~\eqref{edefaction} defines a continuous left action of~$G_\Omega$ on $A_\Omega$.
\end{theorem}

\bp This follows from the proof of \cite[Theorem~6.7]{V}. We include a complete argument for the reader's convenience. Using the identity $(V_\Omega)_{23}(\hat W_\Omega)_{12}(V_\Omega)^*_{23}=(\hat W_\Omega)_{13}(\hat W_\Omega)_{12}$, for $a\in A$ we compute:
\begin{align*}
(\iota\otimes\alpha_\Omega)\eta_\Omega(\alpha(a))
&=(\hat W_\Omega)_{32}(\hat W_\Omega\Omega)_{31}\alpha(a)_{34}(\hat W_\Omega\Omega)^*_{31}(\hat W_\Omega)^*_{32}\\
&=(V_\Omega)_{12}(\hat W_\Omega)_{31}(V_\Omega)^*_{12}\Omega_{31}\alpha(a)_{34}\Omega_{31}^*(V_\Omega)_{12}(\hat W_\Omega)^*_{31}(V_\Omega)^*_{12},
\end{align*}
where $\eta_\Omega$ is defined by \eqref{eetaomega}. Since $(V_\Omega)^*_{12}$ and $\Omega_{31}$ commute, we thus get
$$
(\iota\otimes\alpha_\Omega)\eta_\Omega(\alpha(a))=(V_\Omega)_{12}\eta_\Omega(\alpha(a))_{134}(V_\Omega)_{12}^*.
$$
Multiplying this identity on the left by $1\otimes C_0(G_\Omega)\otimes1\otimes1$, applying the slice maps to the first leg and using that $V_\Omega\in M(K\otimes C_0(G_\Omega))$ and $[(1\otimes C_0(G_\Omega))V^*_\Omega(K\otimes1)]=K\otimes C_0(G_\Omega)$ by regularity of~$G_\Omega$, we see that
$$
[(C_0(G_\Omega)\otimes1\otimes1)\alpha_\Omega\big((T_\nu\otimes\iota)\alpha(A)\big)\mid \nu\in K^*]=C_0(G_\Omega)\otimes [(T_\nu\otimes\iota)\alpha(A)\mid \nu\in K^*].
$$
This implies that $\alpha_\Omega(A_\Omega)\subset M(C_0(G_\Omega)\otimes A_\Omega)$ and that the cancellation property holds. Finally, it is clear that $(\iota\otimes\alpha_\Omega)\alpha_\Omega=(\Delta_\Omega\otimes\iota)\alpha_\Omega$.
\ep

We now want to give a different picture of $A_\Omega$ and $\alpha_\Omega$ based on crossed products. As was mentioned in the introduction, it is inspired by work of Kasprzak~\cite{Kas}, and it was our original motivation for the definition of $A_\Omega$.

Consider the crossed product $G\ltimes_\alpha A$ and the dual action $\hat\alpha$ of $\Ghop$ on $G\ltimes_\alpha A$. We can then try to define a new deformed action $\hat\alpha_\Omega$ of $\Ghopo$ on $G\ltimes_\alpha A$ by
$$
\hat\alpha_\Omega(x)=\Omega_{21}\hat\alpha(x)\Omega^*_{21}.
$$
If $\hat\alpha_\Omega$ is well-defined, $G_\Omega$ is regular and $W_\Omega\otimes1\in M(C_0(G_\Omega)\otimes (G\ltimes_\alpha A))$,
then by a result of Vaes~\cite[Theorem~6.7]{V} discussed in Section~\ref{ssCP}, the action $\hat\alpha_\Omega$ is dual to an action of $G_\Omega$ on a C$^*$-subalgebra $B\subset M(G\ltimes_\alpha A)$, which can be recovered using the homomorphism
\begin{equation}\label{eetaomega1}
\eta_\Omega\colon G\ltimes_\alpha A\to M(K\otimes K\otimes A),\ \ \eta_\Omega(x)=(W_\Omega)^*_{12}\hat\alpha_\Omega(x)(W_\Omega)_{12}.
\end{equation}
Note that since $\hat\alpha(\alpha(a))=1\otimes\alpha(a)$, the definition of $\eta_\Omega$ is consistent with~\eqref{eetaomega}. Note also that for $x=y\otimes1\in C_0(\hat G)\otimes1\in M(G\ltimes_\alpha A)$ we have
$$
\eta_\Omega(x)=(\hat W_\Omega)_{21}(\Dhop_\Omega(y)\otimes1)(\hat W_\Omega)^*_{21}=y\otimes1\otimes1.
$$
It follows that
$$
B=[\{(\omega\otimes\iota\otimes\iota)\eta_\Omega(G\ltimes_\alpha A)\mid \omega\in K^*\}]
=[\{(T_\nu\otimes\iota)\alpha(A)\mid \nu\in K^*\}]=A_\Omega.
$$

Summarizing the above discussion, we have the following result.

\begin{theorem}\label{tVa}
Under the assumptions of Theorem~\ref{tdefaction} suppose that the formula $$\hat\alpha_\Omega(x)=\Omega_{21}\hat\alpha(x)\Omega^*_{21}$$ defines a continuous left action of $\Ghopo$ on $G\ltimes_\alpha A\subset M(K\otimes A)$ and that
$$
W_\Omega\otimes1\in M(C_0(G_\Omega)\otimes(G\ltimes_\alpha A))\subset M(C_0(G_\Omega)\otimes K\otimes A).
$$
Then $G\ltimes_\alpha A=[(C_0(\hat G_\Omega)\otimes1)A_\Omega]$ and the map $\eta_\Omega$ defines an isomorphism $G\ltimes_\alpha A\cong G_\Omega\ltimes_{\alpha_\Omega}A_\Omega$. Under this isomorphism the deformed dual action $\hat\alpha_\Omega$ on $G\ltimes_\alpha A$ becomes the action dual to $\alpha_\Omega$.
\end{theorem}

Again, it is not clear to us what the optimal assumptions for the above two theorems are. Note, however, that for nonregular quantum groups it is not even obvious what the correct definition of a continuous action should be, see the discussion in \cite{BSV}. Even if $G_\Omega$ is regular, it is doubtful that the map $\eta_\Omega$ defines an isomorphism $G\ltimes_\alpha A\cong G_\Omega\ltimes_{\alpha_\Omega}A_\Omega$ for any $A$, since this would imply that the deformed dual action $\hat\alpha_\Omega$ is well-defined on $G\ltimes_\alpha A$, which seems to be overly optimistic already when $\hat G$ is a group and $\Omega$ is a measurable, but not continuous, cocycle on $\hat G$.

\subsection{Deformation in stages}
If $\Omega$ is a cocycle on $\hat G$ and $\Omega_1$ is a cocycle on $\hat G_\Omega$, then it is easy to check that $\Omega_1\Omega$ is a cocycle on $\hat G$. Therefore if the deformed action $\alpha_\Omega$ of $G_\Omega$ on $A_\Omega$ is well-defined, then we can compare the $\Omega_1$-deformation of $A_\Omega$ with the $\Omega_1\Omega$-deformation of $A$.

\begin{theorem}
Assume $G$ is a locally compact quantum group, $\Omega$ is a measurable unitary $2$-cocycle on $\hat G$, $\Omega_1$ is a measurable unitary $2$-cocycle on $\hat G_\Omega$, and $A$ is a C$^*$-algebra equipped with a continuous left action $\alpha$ of $G$. Suppose the following conditions are satisfied:
\enu{i} $A_\Omega=[(T_\nu\otimes\iota)\alpha(A)\mid\nu\in K^*]$;
\enu{ii} the deformed action $\alpha_\Omega$ of $G_\Omega$ is well-defined on $A_\Omega$.

Then the map $x\mapsto(\hat W_\Omega\Omega^*_1)_{21}(1\otimes x)(\hat W_\Omega\Omega^*_1)^*_{21}$ defines an isomorphism $A_{\Omega_1\Omega}\cong(A_\Omega)_{\Omega_1}$. Furthermore, if one of the deformed actions $\alpha_{\Omega_1\Omega}$ and $(\alpha_\Omega)_{\Omega_1}$ is well-defined, then the other is also well-defined and the isomorphism $A_{\Omega_1\Omega}\cong(A_\Omega)_{\Omega_1}$ is $G_{\Omega_1\Omega}$-equivariant.
\end{theorem}

\bp For the proof we need the following identity:
\begin{equation}\label{ecocycle5}
(\hat W_\Omega\Omega)_{23}(\hat W_{\Omega_1\Omega}\Omega_1\Omega)_{12}(\hat W_\Omega\Omega)^*_{23}=(\hat W_{\Omega_1\Omega}\Omega_1)_{12}(\hat W_{\Omega}\Omega)_{13}.
\end{equation}
In order to show this, similarly to the proof of \eqref{ecocycle4} we start with the identity
$$
(\hat W_\Omega\Omega)_{23}(\hat W_\Omega\Omega)_{12}(\hat W_\Omega\Omega)^*_{23}=(\hat W_\Omega)_{12}(\hat W_\Omega\Omega)_{13}.
$$
Next, substitute $(\hat W_\Omega)_{12}$ on both sides for its expression involving $(\hat W_{\Omega_1\Omega})_{12}$ obtained from \eqref{eDC1}:
$$
\hat W_{\Omega_1\Omega}=(\tilde J_1\otimes\hat J)\Omega_1(J_\Omega\otimes\hat J)\hat W_\Omega\Omega_1^*,
$$
where $\tilde J_1$ is the modular involution defined by the dual weight on $W^*(\hat G_\Omega;\Omega_1)$ as explained in Section~\ref{ssDG}. We get
\begin{multline*}
(\hat W_\Omega\Omega)_{23}\big((\tilde J_1\otimes\hat J)\Omega_1(J_\Omega\otimes\hat J)\big)^*_{12}(\hat W_{\Omega_1\Omega}\Omega_1\Omega)_{12}(\hat W_\Omega\Omega)^*_{23}\\=\big((\tilde J_1\otimes\hat J)\Omega_1(J_\Omega\otimes\hat J)\big)^*_{12}((\hat W_{\Omega_1\Omega}\Omega_1)_{12}(\hat W_\Omega\Omega)_{13}.
\end{multline*}
Since $(\hat W_\Omega\Omega)_{23}$ and $\big((\tilde J_1\otimes\hat J)\Omega_1(J_\Omega\otimes\hat J)\big)^*_{12}$ commute, this is exactly \eqref{ecocycle5}.

For $a\in A$ we now start computing:
\begin{equation*}
(\iota\otimes\eta_{\Omega_1}\alpha_\Omega)\eta_\Omega(\alpha(a))=(\hat W_{\Omega_1\Omega}\Omega_1)_{32}(\hat W_\Omega)_{43}\eta_\Omega(\alpha(a))_{145}(\hat W_\Omega)^*_{43}(\hat W_{\Omega_1\Omega}\Omega_1)^*_{32}.
\end{equation*}
By identity \eqref{ecocycle4}, applied to the quantum group $G_\Omega$ and the dual cocycle $\Omega_1$, we have
$$
(\hat W_{\Omega_1\Omega}\Omega_1)_{32}(\hat W_\Omega)_{43}(\hat W_{\Omega_1\Omega}\Omega_1)^*_{32}=(\hat W_\Omega\Omega^*_1)_{43}(\hat W_{\Omega_1\Omega}\Omega_1)_{42}.
$$
Therefore
\begin{align*}
(\iota\otimes\eta_{\Omega_1}\alpha_\Omega)\eta_\Omega(\alpha(a))
&=(\hat W_\Omega\Omega^*_1)_{43}(\hat W_{\Omega_1\Omega}\Omega_1)_{42}\eta_\Omega(\alpha(a))_{145}(\hat W_{\Omega_1\Omega}\Omega_1)^*_{42}(\hat W_\Omega\Omega^*_1)^*_{43}\\
&=(\hat W_\Omega\Omega^*_1)_{43}(\hat W_{\Omega_1\Omega}\Omega_1)_{42}(\hat W_\Omega\Omega)_{41}\alpha(a)_{45}(\hat W_\Omega\Omega)_{41}^*(\hat W_{\Omega_1\Omega}\Omega_1)^*_{42}(\hat W_\Omega\Omega^*_1)^*_{43}.
\end{align*}
By \eqref{ecocycle5} the last expression equals
\smallskip

$\displaystyle
(\hat W_\Omega\Omega^*_1)_{43}(\hat W_\Omega\Omega)_{21}(\hat W_{\Omega_1\Omega}\Omega_1\Omega)_{42}(\hat W_\Omega\Omega)^*_{21}\alpha(a)_{45}(\hat W_\Omega\Omega)_{21}(\hat W_{\Omega_1\Omega}\Omega_1\Omega)^*_{42}(\hat W_\Omega\Omega)^*_{21}(\hat W_\Omega\Omega^*_1)^*_{43}
$
\begin{align*}
&=(\hat W_\Omega\Omega^*_1)_{43}(\hat W_\Omega\Omega)_{21}(\hat W_{\Omega_1\Omega}\Omega_1\Omega)_{42}\alpha(a)_{45}(\hat W_{\Omega_1\Omega}\Omega_1\Omega)^*_{42}(\hat W_\Omega\Omega)^*_{21}(\hat W_\Omega\Omega^*_1)^*_{43}\\
&=(\hat W_\Omega\Omega^*_1)_{43}(\hat W_\Omega\Omega)_{21}\eta_{\Omega_1\Omega}(\alpha(a))_{245}(\hat W_\Omega\Omega)^*_{21}(\hat W_\Omega\Omega^*_1)^*_{43}.
\end{align*}
Thus
$$
(\iota\otimes\eta_{\Omega_1}\alpha_\Omega)\eta_\Omega(\alpha(a))
=(\hat W_\Omega\Omega)_{21}(\hat W_\Omega\Omega^*_1)_{43}\eta_{\Omega_1\Omega}(\alpha(a))_{245}(\hat W_\Omega\Omega^*_1)^*_{43}(\hat W_\Omega\Omega)^*_{21}.
$$
Applying the slice maps to the first two legs we get the first statement of the theorem.

In order to show that the isomorphism $A_{\Omega_1\Omega}\cong(A_\Omega)_{\Omega_1}$ is $G_{\Omega_1\Omega}$-equivariant we need the identity
\begin{equation*}\label{epent}
(\hat W_{\Omega_1\Omega})_{23}(\hat W_\Omega\Omega^*_1)_{12}=(\hat W_\Omega\Omega^*_1)_{12}(\hat W_{\Omega_1\Omega})_{13}(\hat W_{\Omega_1\Omega})_{23}.
\end{equation*}
Since $
\hat W_\Omega\Omega_1^*=\big((\tilde J_1\otimes\hat J)\Omega_1(J_\Omega\otimes\hat J)\big)^*\hat W_{\Omega_1\Omega},
$ this is simply the pentagon relation for $\hat W_{\Omega_1\Omega}$. Denoting the isomorphism in the formulation of the theorem by $\theta$ we compute:
\begin{align*}
(\alpha_{\Omega})_{\Omega_1}(\theta(x))&=(\hat W_{\Omega_1\Omega})_{21}(1\otimes\theta(x))(\hat W_{\Omega_1\Omega})^*_{21}\\
&=(\hat W_{\Omega_1\Omega})_{21}(\hat W_\Omega\Omega^*_1)_{32}(1\otimes1\otimes x)(\hat W_\Omega\Omega^*_1)_{32}^*(\hat W_{\Omega_1\Omega})^*_{21}\\
&=(\hat W_\Omega\Omega^*_1)_{32}(\hat W_{\Omega_1\Omega})_{31}(\hat W_{\Omega_1\Omega})_{21}(1\otimes1\otimes x)(\hat W_{\Omega_1\Omega})_{21}^*(\hat W_{\Omega_1\Omega})_{31}^*(\hat W_\Omega\Omega^*_1)_{32}^*\\
&=(\hat W_\Omega\Omega^*_1)_{32}\alpha_{\Omega_1\Omega}(x)_{13}(\hat W_\Omega\Omega^*_1)_{32}^*=(\iota\otimes\theta)\alpha_{\Omega_1\Omega}(x).
\end{align*}
This proves the second statement of the theorem.
\ep

\begin{example}\mbox{\ }
\enu{i} It is straightforward to check that the deformation $A_1$ of $A$ with respect to the trivial cocycle $1$ is $\alpha(A)$. Therefore given a cocycle $\Omega$ such that the assumptions (i) and (ii) in the above theorem are satisfied, it follows that the map $x\mapsto(\hat W_\Omega\Omega)_{21}(1\otimes x)(\hat W_\Omega\Omega)^*_{21}$ defines an isomorphism $\alpha(A)\cong(A_\Omega)_{\Omega^*}$. In other words,
the map $\eta_\Omega\alpha$ is an isomorphism $A\cong (A_\Omega)_{\Omega^*}$.
\enu{ii} Assume $A=\Ghop\ltimes_\gamma B$ and $\alpha=\hat\gamma$. Then by Proposition~\ref{pdefvstw} we have $A_\Omega\cong \Ghop\ltimes_{\gamma,\Omega}B$. By Proposition~\ref{pdeftw} we have a dual action on $\Ghop\ltimes_{\gamma,\Omega}B$. It is easy to check that this is exactly the deformed action $\alpha_\Omega$. Furthermore, as we have already remarked, the proof of Proposition~\ref{pdefvstw} shows that $A_\Omega=[(T_\nu\otimes\iota)\alpha(A)\mid\nu\in K^*]$. Therefore for dual actions conditions (i) and (ii) in the above theorem are always satisfied. For any cocycle $\Omega_1$ on $\hat G_\Omega$ we thus get
$$
(\Ghop\ltimes_{\gamma,\Omega}B)_{\Omega_1}\cong (\Ghop\ltimes_{\gamma}B)_{\Omega_1\Omega}\cong \Ghop\ltimes_{\gamma,\Omega_1\Omega}B.
$$
In particular, for the C$^*$-algebra $C^*_r(\hat G;\Omega)$ equipped with the action $\beta_\Omega$ of $G_\Omega$ we get $C^*_r(\hat G;\Omega)_{\Omega_1}\cong C^*_r(\hat G;\Omega_1\Omega)$.
\enu{iii} As a particular case of either of the previous two examples we get an isomorphism
$$
C_0(G_\Omega)\cong C^*_r(\hat G_\Omega;\Omega^*)_\Omega,
$$
where the deformation of $C^*_r(\hat G_\Omega;\Omega^*)$ is defined using the action $\beta_{\Omega^*}(x)=W^*(1\otimes x)W$ of~$G$ on $C^*_r(\hat G_\Omega;\Omega^*)$. Explicitly, by the first example the isomorphism $C_0(G_\Omega)\cong (C_0(G_\Omega)_{\Omega^*})_\Omega$ is given by $$x\mapsto \eta_{\Omega^*}\Delta_{\Omega}(x)=(\hat W\Omega^*)_{21}(V_\Omega)_{23}(1\otimes x\otimes 1)(V_\Omega)^*_{23}(\hat W\Omega^*)_{21}^*,$$
so using that $C_0(G_\Omega)_{\Omega^*}=V_\Omega(C^*_r(\hat G_\Omega;\Omega^*)\otimes 1)V_\Omega^*$ we conclude that the isomorphism $C_0(G_\Omega)\cong C^*_r(\hat G_\Omega;\Omega^*)_\Omega$ is given by
$$
x\mapsto(\hat W\Omega^*)_{21}(1\otimes x)(\hat W\Omega^*)_{21}^*.
$$
This is also not difficult to check directly from the definition of $C^*_r(\hat G_\Omega;\Omega^*)_\Omega$.

Note that $C^*_r(\hat G_\Omega;\Omega^*)=\hat J C^*_r(\hat G;\Omega)\hat J$ is nothing other than the deformation of $C_0(G)$ with respect to the right action $\Delta$ of $G$ on $C_0(G)$. More precisely, a right action of $G$ can be considered as a left action of $\Gop$, which is the quantum group $(L^
\infty(G),\Dop)$. The element $(\hat J\otimes\hat J)\Omega(\hat J\otimes\hat J)$ is a dual cocycle on $\Gop$. It is not difficult to check then that the deformation of $C_0(G)$ with respect to the left action of $\Gop$ and the cocycle $(\hat J\otimes\hat J)\Omega(\hat J\otimes\hat J)$ is isomorphic to $C^*_r(\hat G_\Omega;\Omega^*)$, and under this isomorphism the action $\beta_{\Omega^*}$ of $G$ on $C^*_r(\hat G_\Omega;\Omega^*)$ corresponds to the action arising from the left action~$\Delta$ of~$G$ on~$C_0(G)$. Therefore the isomorphism $C_0(G_\Omega)\cong C^*_r(\hat G_\Omega;\Omega^*)_\Omega$ is consistent with what we should expect from the case of finite quantum groups, when $C_0(G_\Omega)\cong C_0(G)$ as coalgebras, while the new algebra structure is obtained by duality from $\Dhat_\Omega=\Omega\Dhat(\cdot)\Omega^{-1}$, which implies that it is obtained by deforming the original product structure on $C_0(G)$ twice, with respect to the left and right actions of $G$ on~$C_0(G)$.
\ee
\end{example}

We finish our general discussion of $\Omega$-deformations with the observation that up to isomorphism the C$^*$-algebra $A_\Omega$ depends only on the cohomology class of $\Omega$. Recall that given a $2$-cocycle $\Omega$ on~$\hat G$ and a unitary $u\in L^\infty(\hat G)$, the element $\Omega_u=(u\otimes u)\Omega\Dhat(u)^*$ is again a cocycle on $\hat G$. The cocycles~$\Omega$ and~$\Omega_u$ are called cohomologous. The set of cohomology classes of unitary $2$-cocycles on~$\hat G$ is denoted by $H^2(\hat G;\T)$. In general it is just a set.

\begin{proposition} \label{pCob}
Assume $G$ is a locally compact quantum group, $\Omega$ is a measurable unitary $2$-cocycle on $\hat G$ and $u$ is a unitary in $L^\infty(\hat G)$. Then for any C$^*$-algebra $A$ equipped with a continuous left action of $G$, the map $\Ad(u\otimes 1)$ defines an isomorphism $A_\Omega\cong A_{\Omega_u}$.
\end{proposition}

\bp Since $\hat W\Omega_u^*=\hat W\Dhat(u)\Omega^*(u^*\otimes u^*)=(1\otimes u)\hat W\Omega^*(u^*\otimes u^*)$, the map $\Ad u$ defines a right $G$-equivariant isomorphism $W^*(\hat G;\Omega)\cong W^*(\hat G;\Omega_u)$. This isomorphism maps $(\omega(\cdot\, u^*)\otimes\iota)(\hat W\Omega^*)$ into $(\omega\otimes\iota)(\hat W\Omega^*_u)$. Consider the GNS-representations $\tilde\Lambda\colon\mathcal N_{\tilde\varphi}\to L^2(G)$ and $\tilde\Lambda_u\colon\mathcal N_{\tilde\varphi_u}\to L^2(G)$ defined by the dual weights $\tilde\varphi$ on $W^*(\hat G;\Omega)$ and $\tilde\varphi_u$ on $W^*(\hat G;\Omega_u)$, as described in Section~\ref{ssDG}. Then the isomorphism $W^*(\hat G;\Omega)\cong W^*(\hat G;\Omega_u)$ defines a unitary $\tilde u$ on $L^2(G)$ such that $\tilde u\tilde\Lambda(x)=\tilde\Lambda_u(uxu^*)$ for $x\in\mathcal N_{\tilde\varphi}$,~so
$$
\tilde u\tilde\Lambda((\omega(\cdot\, u^*)\otimes\iota)(\hat W\Omega^*))=
\tilde\Lambda_u((\omega\otimes\iota)(\hat W\Omega^*_u))
$$
for suitable $\omega\in K^*$. Since $\tilde\Lambda((\omega(\cdot\, u^*)\otimes\iota)(\hat W\Omega^*))=\Lambda((\omega(\cdot\, u^*)\otimes\iota)(\hat W))$ and a similar formula holds for $\tilde\Lambda_u$, we in other words have
$$
\tilde u\Lambda((\omega(\cdot\, u^*)\otimes\iota)(\hat W))=\Lambda(\omega\otimes\iota)(\hat W)).
$$
But $\Lambda((\omega(\cdot\, u^*)\otimes\iota)(\hat W))=u^*\Lambda((\omega\otimes\iota)(\hat W))$, as can be easily checked using that $$(\Lambda((\nu\otimes\iota)(\hat W)),\hat\Lambda(y))=\nu(y^*).$$ It follows that $\tilde u=u$. Hence the modular involution $\tilde J_u$ defined by the  weight $\tilde\varphi_u$ on $W^*(\hat G;\Omega_u)$ is equal to $u\tilde Ju^*$. Therefore
$$
(\tilde J_u\otimes\hat J)\Omega_u(J\otimes\hat J)
=(u\otimes\hat Ju\hat J)\big((\tilde J\otimes\hat J)\Omega(J\otimes\hat J)\big)(J\otimes\hat J)\Dhat(u)^*(J\otimes\hat J).
$$
Next, we have
$$
(J\otimes\hat J)\Dhat(u)^*(J\otimes\hat J)\hat W
=(J\otimes\hat J)\Dhat(u)^*\hat W^*(J\otimes\hat J)
=(J\otimes\hat J)\hat W^*(1\otimes u^*)(J\otimes\hat J)
=\hat W(1\otimes\hat Ju^*\hat J).
$$
Hence, by \eqref{eDC1},
$$
\hat W_{\Omega_u}\Omega_u=(\tilde J_u\otimes\hat J)\Omega_u(J\otimes\hat J)\hat W=(u\otimes \hat Ju\hat J)\hat W_\Omega\Omega(1\otimes\hat Ju^*\hat J).
$$
Recalling the definition of $\eta_\Omega$ we get
$$
\eta_{\Omega_u}=\Ad(\hat Ju\hat J\otimes u\otimes1)\eta_\Omega.
$$
This gives the result.
\ep

\bigskip

\section{Cocycles on group duals} \label{sGD}

In this section we assume that $G$ is a genuine locally compact group.

\subsection{Dual cocycles and deformations of the Fourier algebra}\label{ssFourier}
Denote by $\lambda_g$, resp. $\rho_g$, the operators of the left, resp. right, regular representation of $G$. Then
$$
(\hat W\xi)(s,t)=\xi(ts,t)=(\lambda_t^{-1}\xi(\cdot,t))(s)\ \ \text{and}\ \ (V\xi)(s,t)=(\rho_t\xi(\cdot,t))(s)\ \ \text{for}\ \ \xi\in L^2(G\times G).
$$
The predual of $L^\infty(\hat G)=W^*(G)$ can be identified with the Fourier algebra $A(G)\subset C_0(G)$, so an element $\omega\in W^*(G)_*$ is identified with the function $f(g)=\omega(\lambda_g)$ on $G$. Note that under this identification we have
$$
(f\otimes\iota)(\hat W)=\check{f}\ \ \text{for}\ \ f\in A(G),
$$
where $\check{f}(g)=f(g^{-1})$.

Assume now that $\Omega$ is a measurable unitary $2$-cocycle on $\hat G$. Then we can define a new product~$\star_\Omega$ on $A(G)$ by
$$
f_1\star_\Omega f_2=(f_1\otimes f_2)(\Dhat(\cdot)\Omega^*).
$$
The associativity of this product is equivalent to the cocycle identity for~$\Omega$. Identity~\eqref{ecocycle0} shows that the formula
$$
\pi_\Omega(f)=(f\otimes\iota)(\hat W\Omega^*)
$$
defines a representation of $(A(G),\star_\Omega)$ on $L^2(G)$, and then by definition the C$^*$-algebra $C^*_r(\hat G;\Omega)$ is generated by $\pi_\Omega(A(G))$. Recall that by Theorem~\ref{tCdef} we in fact have $C^*_r(\hat G;\Omega)=\overline{\pi_\Omega(A(G))}$. Nevertheless we do not claim that $(A(G),\star_\Omega)$ is itself a $*$-algebra, although this is often the case.  In Section~\ref{sBG} we will give an example where $\pi_\Omega(A(G))$ is not a $*$-subalgebra of $C^*_r(\hat G;\Omega)$. Since by~\eqref{etildelambda}, we have
$$
\tilde\Lambda((f\otimes\iota)(\hat W\Omega^*))=\Lambda((f\otimes\iota)(\hat W))=\check f,
$$
the representation $\pi_\Omega$ is given by
$$
\pi_\Omega(f_1)\check f_2=(f_1\star_\Omega f_2)\check{}\ \ \text{for}\ \ f_1\in A(G)\ \ \text{and}\ \ f_2\in A(G)\cap C_c(G).
$$
In other words, since $\int_G\check f(g)dg=\int_Gf(g)\Delta_G(g)^{-1}dg$, the representation $\pi_\Omega$ is simply the left regular representation of $(A(G),\star_\Omega)$ on itself, with $A(G)$, or more precisely $A(G)\cap C_c(G)$, completed to a Hilbert space using the scalar product defined by the right Haar measure $\Delta_G(g)^{-1}dg$.

The left translations of $G$ on itself define automorphisms of $(A(G),\star_\Omega)$. On the level of $C^*_r(\hat G;\Omega)$ the action by left translations is exactly the action $\beta$ introduced earlier, so
$$
\beta_g(\pi_\Omega(f))=\pi_\Omega(f(g^{-1}\cdot))\ \ \text{for}\ \ f\in A(G).
$$
This is the reason for the appearance of the right Haar measure, since the average of a function on~$G$ with respect to the action by left translations is the integral with respect to a right Haar measure.

\smallskip

Conversely, assume we have a product $\star$ on $A(G)$ that is invariant under left translations. Assume there exists an element $\Omega\in W^*(G)\bar\otimes W^*(G)$ such that
$$
(f_1\otimes f_2)(\Omega^*)=(f_1\star f_2)(e)\ \ \text{for all}\ \ f_1,f_2\in A(G),
$$
which happens exactly when the map $f_1\otimes f_2\mapsto (f_1\star f_2)(e)$ extends to a bounded linear functional on the projective tensor product $A(G)\hat\otimes A(G)$. Then $f_1\star f_2=(f_1\otimes f_2)(\Dhat(\cdot)\Omega^*)$. For finite groups this was observed by Movshev~\cite{Mo}. As follows from results in \cite{Mo}, if $G$ is finite and $(A(G),\star)$ is semisimple, then $\Omega$ is invertible and cohomological to a unitary cocycle, that is, there exists an invertible element $a\in W^*(G)$ such that $(a\otimes a)\Omega\Dhat(a)^{-1}$ is unitary.

\smallskip

For a related discussion see~\cite{LR1}.

\subsection{$\Omega$-Deformations and generalized fixed point algebras}\label{ssFixed}
As we know, by \cite[Lemma~1.12]{VV} the action $\beta$ of $G$ on $W^*(\hat G;\Omega)$ is integrable. We will now show a stronger property: the action of~$G$ on $C^*_r(\hat G;\Omega)$ is integrable. Since the integrability property that we will establish appears under several different names in the literature, let us say precisely what we mean by this.

Assume $\gamma$ is a continuous action of $G$ on a C$^*$-algebra $B$. An element $b\in B_+$ is called $\gamma$-integrable if there exists an element $\psi_\gamma(b)\in M(B)$ such that for every state $\omega$ on $B$ the function $g\mapsto\omega(\gamma_g(b))$ is integrable and
$$
\int_G\omega(\gamma_g(b))dg=\omega(\psi_\gamma(b)).
$$
Clearly, the element $\psi_\gamma(b)$ is uniquely determined by $b$, and $\psi_\gamma(b)\in M(B)^\gamma$.

There are several other equivalent definitions of integrable elements, see~\cite{Ri98} (note that in \cite{Ri98} $\gamma$-integrable elements are called $\gamma$-proper). For example, by \cite[Proposition~4.4]{Ri98} a positive element $b$ is $\gamma$-integrable if and only if the functions $g\mapsto\gamma_g(b)c$ and $g\mapsto c\gamma_g(b)$ are unconditionally integrable for all $c\in B$, meaning that their integrals over compact subsets of $G$ form Cauchy nets.

The set $\PP^+_\gamma$ of $\gamma$-integrable positive  elements is a hereditary cone in $B_+$, see e.g.~\cite[Lemma~2.7]{Ri98}. Hence the linear span of $\PP^+_\gamma$ is a $*$-algebra. We say that $\gamma$ is integrable if $\PP^+_\gamma$ is dense in $B_+$, or equivalently, $[\PP^+_\gamma]=B$.

\begin{proposition} \label{pintegr}
For any measurable unitary $2$-cocycle $\Omega$ on $\hat G$ the action $\beta$ of $G$ on $C^*_r(\hat G;\Omega)$ is integrable.
\end{proposition}

Before we turn to the proof, let us discuss the difference between this statement and the integrability of the action of $G$ on $W^*(\hat G;\Omega)$. By \cite[Theorem~1.11]{VV} the action on $W^*(\hat G;\Omega)$ is ergodic. In particular, $M(C^*_r(\hat G;\Omega))^\beta=\C1$. For $x\in W^*(\hat G;\Omega)_+$  in the domain of definition $\operatorname{dom}\tilde\varphi$ of the dual weight $\tilde\varphi$ we have
\begin{equation}\label{eintegr}
\int_G\omega(\beta_g(x))dg=\tilde\varphi(x)
\end{equation}
for all normal states $\omega$ on $W^*(\hat G;\Omega)$. It follows that if $x\in C^*_r(\hat G;\Omega)_+$ is $\beta$-integrable, then $x\in \operatorname{dom}\tilde\varphi$ and $\psi_\beta(x)=\tilde\varphi(x)1$. The difference between the cones $\PP^+_\beta$ and $C^*_r(\hat G;\Omega)\cap\operatorname{dom}\tilde\varphi$ is that for the elements of $\PP^+_\beta$ identity \eqref{eintegr} should be satisfied for all states $\omega$ on $C^*_r(\hat G;\Omega)$, while for the elements of $C^*_r(\hat G;\Omega)\cap\operatorname{dom}\tilde\varphi$ we need only to consider normal states on $W^*(\hat G;\Omega)$. Note also that the density of $C^*_r(\hat G;\Omega)\cap\operatorname{dom}\tilde\varphi$ in $C^*_r(\hat G;\Omega)_+$ follows already from the proof of \cite[Lemma~1.12]{VV}.

\bp[Proof of Proposition~\ref{pintegr}]
Let $\nu$ be a normal state on $W^*(\hat G_\Omega;\Omega^*)$.
As we already observed in Section~\ref{squantmaps}, the quantization map $T_\nu\colon M(C_0(G))\to M(C^*_r(\hat G;\Omega))$ is strictly continuous on bounded sets. It follows that for any state $\omega$ on $C^*_r(\hat G;\Omega)$, the positive linear functional $\omega T_\nu$ on $C_0(G)$ is again a state.

Since $T_\nu\colon L^\infty(G)\to W^*(\hat G;\Omega)$ is a $G$-equivariant normal u.c.p.~map, if $f\in\operatorname{dom}\varphi=L^\infty(G)_+\cap L^1(G)$, then $T_\nu(f)\in \operatorname{dom}\tilde\varphi$ and
$$
\tilde\varphi(T_\nu(f))=\varphi(f)=\int_Gf(g)dg.
$$
For the action of $G$ on itself by right translations there is no distinction between integrability of an element $f\in C_0(G)_+$ in the von Neumann algebraic and the C$^*$-algebraic sense: both conditions are equivalent to $f\in L^1(G)$. For $f\in C_0(G)_+\cap L^1(G)$ and any state $\omega$ on $C^*_r(\hat G;\Omega)$ we then get
$$
\int_G\omega(\beta_g(T_\nu(f)))dg=\int_G\omega T_\nu(f(\cdot\,g))dg=(\omega T_\nu)(1)\varphi(f)=\tilde\varphi(T_\nu(f)),
$$
so $T_\nu(f)$ is $\beta$-integrable.

By Proposition~\ref{pquantmaps} the span of the spaces $T_\nu(C_0(G))$, $\nu\in W^*(\hat G_\Omega;\Omega^*)_*$, is dense in  $C^*_r(\hat G;\Omega)$. Hence $[\PP^+_\beta]=C^*_r(\hat G;\Omega)$.
\ep

Returning to a general action $\gamma$ of $G$ on $B$, it is easy to see that if $b\in B_+$ is $\gamma$-integrable and $x\in M(B)^\gamma$, then $x^*bx$ is again $\gamma$-integrable and $\psi_\gamma(x^*bx)=x^*\psi_\gamma(b)x$. This implies that the span of $\PP^+_\gamma$ is an $M(B)^\gamma$-bimodule, which in turn implies that the span of $\psi_\gamma(\PP^+_\gamma)$ is a $*$-ideal in $M(B)^\gamma$. The C$^*$-algebra $[\psi_\gamma(\PP^+_\gamma)]\subset M(B)^\gamma$ can be considered as a generalized fixed point algebra for the action $\gamma$ on $B$. In general, however, it is too big to have good properties and it is not clear what the correct definition of a generalized fixed point algebra should be, see the discussion in~\cite{Ri98}.

\smallskip

Consider now a continuous action $\alpha$ of $G$ on a C$^*$-algebra $A$. By the observation immediately after Definition~\ref{domegadef}, we have $A_\Omega\subset M(C^*_r(\hat G;\Omega)\otimes A)^{\beta\otimes\alpha}$. We can now prove a slightly more precise result.

\begin{proposition} \label{pfixed2}
For any measurable unitary $2$-cocycle $\Omega$ on $\hat G$ and any C$^*$-algebra $A$ equipped with a continuous action $\alpha$ of $G$, the diagonal action $\beta\otimes\alpha$ of $G$ on $C^*_r(\hat G;\Omega)\otimes A$ is integrable and the C$^*$-algebra $A_\Omega$ is contained in $[\psi_{\beta\otimes\alpha}(\PP^+_{\beta\otimes\alpha})]\subset M(C^*_r(\hat G;\Omega)\otimes A)^{\beta\otimes\alpha}$.
\end{proposition}

\bp The first statement follows immediately from the integrability of $\beta$.

In order to prove the second statement denote by $\rho$ the action of $G$ on $C_0(G)$ by right translations. Then it is easy to check that for any $f\in C_0(G)_+\cap L^1(G)$ and $a\in A_+$ the element $f\otimes a$ is $\rho\otimes\alpha$-integrable and
$$
\psi_{\rho\otimes\alpha}(f\otimes a)=\alpha(a_f),
$$
where $a_f=\int_Gf(g)\alpha_g(a)dg$; note that if we identify $M(C_0(G)\otimes A)$ with $C_b(G;M(A))$ then by definition $\alpha(a_f)(g)=\alpha_{g^{-1}}(a_f)$.

From this, arguing as in the proof of the previous proposition, we conclude that the element $T_\nu(f)\otimes a$ is $\beta\otimes\alpha$-integrable and
$$
\psi_{\beta\otimes\alpha}(T_\nu(f)\otimes a)=(T_\nu\otimes\iota)\alpha(a_f).
$$
Since $[\psi_{\beta\otimes\alpha}(\PP^+_{\beta\otimes\alpha})]$ is a C$^*$-algebra, it follows that $A_\Omega$ is contained in $[\psi_{\beta\otimes\alpha}(\PP^+_{\beta\otimes\alpha})]$.
\ep

In the case when $G$ is a compact group all the analytical difficulties disappear and we get the following.

\begin{proposition}
If $G$ is compact, then $A_\Omega=(C^*_r(\hat G;\Omega)\otimes A)^{\beta\otimes\alpha}$.
\end{proposition}

\bp When $G$ is compact, every positive element of $C^*_r(\hat G;\Omega)\otimes A$ is integrable and the map $\psi_{\beta\otimes\alpha}$ extends by linearity to a bounded map $C^*_r(\hat G;\Omega)\otimes A\to M(C^*_r(\hat G;\Omega)\otimes A)^{\beta\otimes\alpha}$ with image $(C^*_r(\hat G;\Omega)\otimes A)^{\beta\otimes\alpha}$. As follows from the previous proof, a dense subspace of elements of $C^*_r(\hat G;\Omega)\otimes A$ is mapped by $\psi_{\beta\otimes\alpha}$ onto a generating subspace of $A_\Omega$. Hence $A_\Omega=(C^*_r(\hat G;\Omega)\otimes A)^{\beta\otimes\alpha}$.
\ep

For an overview of what is known about cocycles on duals of compact groups see \cite{NTijm}.

\subsection{Regularity of cocycles and proper actions}\label{ssRegular2}
A stronger notion than integrability was introduced by Rieffel in \cite{Ri}. Namely, an action $\gamma$ of $G$ on a C$^*$-algebra $B$ is called proper, if there exists a dense $\gamma$-invariant $*$-subalgebra $B_0\subset B$ such that for all $b,c\in B$ the functions $g\mapsto b\gamma_g(c)$ and $g\mapsto\Delta_G(g)^{1/2}b\gamma_g(c)$ are norm-integrable\footnote{We define the modular function so that $\int_Gf(gh)dg=\Delta_G(h)^{-1}\int_Gf(g)dg$, which is opposite to the conventions in~\cite{Ri}.} and there exists an element $x\in M(B_0)^\gamma\subset M(B)^\gamma$ such that $\int_Ga\gamma_g(bc)dg=ax$ for all $a\in B_0$. By \cite[Proposition~4.6]{Ri98} a proper action is integrable; in fact, if $b\in B_0$ then $b^*b$ is $\gamma$-integrable.

The integrable functions $g\mapsto\Delta_G(g)^{1/2}b\gamma_g(c)$ define elements of the reduced crossed product $G\ltimes_\gamma B$. The closure $I$ of the space spanned by such elements is a $*$-ideal in $G\ltimes_\gamma B$. As shown in~\cite{Ri}, this ideal is strongly Morita equivalent to the C$^*$-subalgebra $[\psi_\gamma(B^2_0)]\subset M(B)^\gamma$. The action~$\gamma$ is called saturated if $B_0$ can be chosen such that $I=G\ltimes_\gamma B$.

The relevance of these notions for us is explained by the following.

\begin{proposition}\label{pregular}
For a measurable unitary $2$-cocycle $\Omega$ on $\hat G$, assume the action $\beta$ of $G$ on $C^*_r(\hat G;\Omega)$ is proper and saturated. Then $\Omega$ is regular.
\end{proposition}

\bp Since $M(C^*_r(\hat G;\Omega))^\beta=\C1$, by the above discussion the assumptions of the proposition imply that $C^*_r(\hat G)\rtimes_\beta G$ is strongly Morita equivalent to $\C$, that is, $C^*_r(\hat G;\Omega)\rtimes_\beta G$ is isomorphic to the algebra of compact operators on some Hilbert space. But this is one of the equivalent characterizations of regularity.
\ep

We expect the saturation property to hold more or less automatically. For example, it holds when~$G$~is compact, in which case, however, we do need the above proposition to show regularity.

\smallskip

We finish this section with a simple result on continuity of dual cocycles.

\begin{proposition}
Assume $\Omega$ is a measurable unitary $2$-cocycle on $\hat G$ such that both $\Omega$ and $\Omega^*$ map $C_c(G\times G)\subset L^2(G\times G)$ into $L^1(G\times G)\cap L^2(G\times G)$. Then $\Omega$ is continuous, that is, $\Omega\in M(C^*_r(G)\otimes C^*_r(G))$.
\end{proposition}

\bp A function $f\in L^1(G\times G)\cap L^2(G\times G)$ defines both a vector in $L^2(G\times G)$ and an element $(\lambda\otimes\lambda)(f)\in C^*_r(G)\otimes C^*_r(G)$. We claim that if $f\in C_c(G\times G)$, then
$$
\Omega(\lambda\otimes\lambda)(f)=(\lambda\otimes\lambda)(\Omega f).
$$
Indeed, if $\xi\in C_c(G\times G)$, then, using that $\Omega$ commutes with the operator $\zeta\mapsto\zeta*\xi$ on $L^2(G\times G)$, we have
$$
\Omega(\lambda\otimes\lambda)(f)\xi=\Omega(f*\xi)=(\Omega f)*\xi=(\lambda\otimes\lambda)(\Omega f)\xi.
$$
Since $C_c(G\times G)$ is dense in $L^2(G\times G)$, this proves our claim.

It follows that $\Omega(C^*_r(G)\otimes C^*_r(G))\subset C^*_r(G)\otimes C^*_r(G)$. Similarly, $\Omega^*(\lambda\otimes\lambda)(f)=(\lambda\otimes\lambda)(\Omega^* f)$ for all $f\in C_c(G\times G)$, which implies that $\Omega^*(C^*_r(G)\otimes C^*_r(G))\subset C^*_r(G)\otimes C^*_r(G)$.
\ep

\bigskip

\section{Dual cocycles for a class of solvable Lie groups} \label{sBG}

In this section we briefly consider dual cocycles recently constructed by Bieliavsky et al.~\cite{BG,Betal}.

\subsection{Deformation of negatively curved K\"ahlerian Lie groups}
As explained in \cite{BG}, by results of Pyatetskii-Shapiro, negatively curved K\"ahlerian Lie groups can be decomposed into iterated semidirect products of certain elementary groups, called elementary normal $j$-groups in~\cite{BG}. To simplify matters we will consider only the latter groups. Thus, throughout the whole Section~\ref{sBG} we assume that $G$ is a simply connected real Lie group of dimension $2d+2$ with a basis $H$, $\{X_j\}^{2d}_{j=1}$, $E$ of the Lie algebra $\g$ satisfying the relations
$$
[H,E]=2E,\ \ [H,X_j]=X_j,\ \ [E,X_j]=0,\ \ [X_i,X_j]=(\delta_{i+d,j}-\delta_{i,j+d})E.
$$
The map
$$
\R\times\R^{2d}\times\R\ni(a,v,t)\mapsto \exp(aH)\exp\left(\sum^{2d}_{j=1}v_jX_j+tE\right)\in G
$$
is a diffeomorphism. In the coordinates $(a,v,t)$ the group law on $G$ takes the form
$$
(a,v,t)(a',v',t')=(a+a',e^{-a'}v+v',e^{-2a'}t+t'+\frac{1}{2}e^{-a'}\omega_0(v,v')),
$$
where $\omega_0(v,v')=\sum^d_{i=1}(v_iv'_{i+d}-v_{i+d}v'_i)$ is the standard symplectic form on $\R^{2d}$. From this formula it is clear that the usual Lebesgue measure on $\R^{2d+2}$ defines a left Haar measure on $G$. Then the modular function $\Delta_G$ is given by\footnote{Note again that our definition of the modular function is opposite to the one in~\cite{BG}.}
$$
\Delta_G(a,v,t)=e^{-(2d+2)a}.
$$

For every $\theta\in\R$, $\theta\ne0$, Bieliavsky and Gayral~\cite[Section~4.1]{BG} construct a new product $\star_\theta$ on a space $\E_\theta(G)$ of smooth functions on $G$. The precise definition of $\E_\theta(G)$ is not important to us. What we need to know is that $\E_\theta(G)$ contains $C^\infty_c(G)$ and there exists a bijective linear map $T_\theta\colon S(\R^{2d+2})\to \E_\theta(G)$ that is compatible with complex conjugation and that extends to a unitary operator $L^2(\R^{2d+2})\to L^2(G)$. The new product is then defined by
$$
f_1\star_\theta f_2=T_\theta(T^{-1}_\theta(f_1)\star^0_\theta T^{-1}_\theta(f_2)),
$$
where $\star^0_\theta$ denotes the standard Moyal product on $S(\R^{2d+2})$ defined using the symplectic form
$$
\omega_\g((a,v,t),(a',v',t'))=2(at'-ta')+\omega_0(v,v').
$$
One of the reasons to introduce the map $T_\theta$ is that the product $\star_\theta$ becomes left $G$-invariant. Furthermore, it is possible to explicitly write down the distribution kernel of the product:
$$
(f_1\star_{\theta} f_2)(g)=\int_{G\times G}K_\theta(x,y)f_1(gx)f_2(gy)dx\,dy\ \ \text{for}\ \ f_1,f_2\in C^\infty_c(G),
$$
where
$$
K_\theta(x,y)=\frac{4}{(\pi\theta)^{2d+2}}A(x,y)\exp\left\{{\frac{2i}{\theta}S(x,y)}\right\}
$$
and, for $x=(a,v,t)$ and $x'=(a',v',t')$,
$$
A(x,x')=\big(\cosh(a)\cosh(a')\cosh(a-a')\big)^d\big(\cosh(2a)\cosh(2a')\cosh(2a-2a')\big)^{1/2},
$$
$$
S(x,x')=\sinh(2a)t'-\sinh(2a')t+\cosh(a)\cosh(a')\omega_0(v,v').
$$

In view of our discussion of the relation between dual cocycles and deformations of the Fourier algebra in Section~\ref{ssFourier} it is then natural to try to define a cocycle $\Omega_\theta$ on $\hat G$ by
$$
\Omega^*_\theta=\int_{G\times G}K_\theta(x,y) \lambda_x\otimes\lambda_y\,dx\,dy\ \ \text{on}\ \ C^\infty_c(G\times G)\subset L^2(G\times G).
$$
The question is whether this defines a unitary operator on $L^2(G\times G)$. By \cite{Betal} this is indeed the case. In fact, the proof is essentially contained already in \cite{BG}. Namely, as established in the proof of \cite[Proposition 8.45]{BG}, we have
$$
\int_{G\times G} K_{\theta} (x,y)K_{-\theta} (gx,hy)\Delta_G(x)\Delta_G(y)dx\,dy=\delta_e (g)\delta_e (h),
$$
with the integral understood in the distribution sense.
It is not difficult to check that this, together with $\overline{K_\theta(x,y)}=K_{-\theta}(x,y)$, implies that $\Omega_\theta$ is an isometry. Similarly, we have
$$
\int_{G\times G} K_{-\theta} (x,y)K_\theta (xg,yh)dx\,dy=\delta_e (g)\delta_e (h),
$$
which implies that $\Omega^*_\theta$ is an isometry.

Therefore, by \cite{Betal}, we get a family of unitary $2$-cocycles $\Omega_\theta$ on $\hat G$. The corresponding quantum groups~$G_{\Omega_\theta}$ provide a nonformal deformation of the Poisson-Lie group $G$, with the Poisson structure defined by the nondegenerate $2$-cocycle $\omega_\g$ on $\g$.

In the simplest case $d=0$ the group $G$ is the $ax+b$ group. A quantization of the same Poisson-Lie structure on it has been defined by Baaj and Skandalis~\cite{Sk}, see also~\cite[Section~5.3]{VV} and~\cite{St} (note also that there exists only one, up to isomorphism and rescaling, different Poisson-Lie structure; it has been quantized by Pusz, Woronowicz and Zakrzewski~\cite{PW,WZ}). It would be interesting to check whether for $d=0$ the quantum groups $G_{\Omega_\theta}$ are isomorphic to the one defined by Baaj and Skandalis.

\subsection{Involution on the twisted group algebra} Given the cocycle $\Omega_\theta$, we can now consider the new product $\star_{\Omega_\theta}$ on the Fourier algebra $A(G)$ and the representation $\pi_{\Omega_\theta}$ of $(A(G),\star_{\Omega_\theta})$ on $L^2(G)$ given by
$$
\pi_{\Omega_\theta}(f_1)\check f_2=(f_1\star_{\Omega_\theta} f_2)\check{}\ \ \text{for}\ \ f_1\in A(G)\ \ \text{and}\ \ f_2\in A(G)\cap C_c(G).
$$
At a first glance a bit surprisingly, the algebra $\pi_{\Omega_\theta}(A(G))$ fails to be a $*$-algebra. Namely, consider the modular function $\Delta_G$ as the unbounded operator of multiplication by $\Delta_G$ on~$L^2(G)$, so identify~$\Delta_G$ with the modular operator defined by the Haar weight on $L^\infty(\hat G)=W^*(G)$. Consider also the dense subspace $A_\infty(G)$ of $A(G)$ spanned by the functions of the form $\xi*\zeta$ with $\xi,\zeta\in C^\infty_c(G)$. We then have the following.

\begin{lemma} \label{lstar}
For any $f\in A_\infty(G)$ we have
$$
\pi_{\Omega_\theta}(f)^*=\Delta_G^{-1}\pi_{\Omega_\theta}(\bar f)\Delta_G\ \ \text{on}\ \ C_c(G)\subset L^2(G).
$$
\end{lemma}

\bp Since for $f\in A_\infty(G)$ we have
$$
(f\otimes\iota)(\hat W(\lambda_x\otimes\lambda_y))=\check f(x^{-1}\cdot)\lambda_y,
$$
we get
$$
\pi_{\Omega_\theta}(f)=\int_{G\times G}K_\theta(x,y)\check f(x^{-1}\cdot)\lambda_y\,dx\,dy\ \ \text{on}\ \ C_c(G).
$$
It follows that on $C_c(G)$ we have
\begin{align*}
\pi_{\Omega_\theta}(f)^*&=\int_{G\times G}\overline{K_\theta(x,y)}\lambda_{y^{-1}}\check{\bar f}(x^{-1}\cdot)dx\,dy
=\int_{G\times G}\overline{K_\theta(x,y)}\check{\bar f}(x^{-1}y\,\cdot)\lambda_{y^{-1}}\,dx\,dy\\
&=\int_{G\times G}\overline{K_\theta(yx,y)}\check{\bar f}(x^{-1}\cdot)\lambda_{y^{-1}}\,dx\,dy
=\int_{G\times G}\overline{K_\theta(y^{-1}x,y^{-1})}\Delta_G(y)^{-1}\check{\bar f}(x^{-1}\cdot)\lambda_y\,dx\,dy\\
&=\Delta_G^{-1} \int_{G\times G}\overline{K_\theta(y^{-1}x,y^{-1})}\check{\bar f}(x^{-1}\cdot)\lambda_y\,dx\,dy\, \Delta_G.
\end{align*}
Therefore it suffices to check that $\overline{K_\theta(y^{-1}x,y^{-1})}=K_\theta(x,y)$, or equivalently,
$$
A(y^{-1}x,y^{-1})=A(x,y)\ \ \text{and}\ \ S(y^{-1}x,y^{-1})=-S(x,y).
$$
Both identities are checked by a straightforward computation.
\ep

We will now give a different proof of this lemma using known properties of the Moyal product. As we discussed in Section~\ref{ssFourier}, the representation $\pi_{\Omega_\theta}$ can be thought of as the left regular representation of $(A(G),\star_{\Omega_\theta})$ with respect to the scalar product defined by the right Haar measure. We can also try to use the left Haar measure. In general we see no reason to expect the corresponding representation to be well-defined. But in the present case, where $\star_\theta$ was constructed using the Moyal product, we do have a representation $\pi_\theta$ of $(\E_\theta(G),\star_\theta)$ on $L^2(G)$ defined by
$$
\pi_\theta(f_1)f_2=f_1\star_\theta f_2\ \ \text{for}\ \ f_1,f_2\in\E_\theta(G).
$$
Furthermore, for this representation we have $\pi_\theta(f)^*=\pi_\theta(\bar f)$. Since by construction the products~$\star_\theta$ and~$\star_{\Omega_\theta}$ coincide on $A_\infty(G)\subset A(G)\cap\E_\theta(G)$ (but $A_\infty(G)$ is not closed under these products), for any $f\in A_\infty(G)$ we have
$$
\pi_{\Omega_\theta}(f)=\nabla\pi_\theta(f)\nabla\ \ \text{on}\ \ A_\infty(G)\subset L^2(G),
$$
where $\nabla$ is the unbounded involutive operator defined by $\nabla f=\check f$. Since $\nabla^*=\Delta_G^{-1}\nabla=\nabla\Delta_G$, this is consistent with Lemma~\ref{lstar}.

In view of the identity $\pi_{\Omega_\theta}(f)=\nabla\pi_\theta(f)\nabla$ it may seem surprising that both representations $\pi_{\Omega_\theta}$ and $\pi_\theta$ are well-defined. The representation $\pi_{\Omega_\theta}$ is well-defined by our general theory. The reason why $\pi_\theta$ is well-defined is that ultimately the product $\star_\theta$ was constructed using a dual cocycle on~$\R^{2d+2}$ and this group is unimodular.

Turning to C$^*$-algebras, the obvious conclusion is that the identity map on $A_\infty(G)$ does not extend to a $*$-isomorphism of $C^*_r(\hat G;\Omega_\theta)=\overline{\pi_{\Omega_\theta}(A(G))}$ and $\overline{\pi_\theta(\E_\theta(G))}$. This, however, does not exclude the possibility that these C$^*$-algebras are isomorphic in a canonical $G$-equivariant way. In fact, the above considerations suggest that the conjugation by the involutive unitary $\Delta_G^{-1/2}\nabla=J\hat J$ gives such an isomorphism. This will be analyzed in a subsequent publication.

The C$^*$-algebra $\overline{\pi_\theta(\E_\theta(G))}$ is the $\theta$-deformation of $C_0(G)$ as defined by Bieliavsky and Gayral~\cite{BG}. To be more precise, instead of the representation $\pi_\theta$ they use the Weyl quantization map. But it is well-known that this gives a quasi-equivalent representation. In particular, as an abstract C$^*$-algebra, $\overline{\pi_\theta(\E_\theta(G))}$ is isomorphic to the algebra of compact operators on an infinite dimensional separable Hilbert space.

\subsection{Two-parameter deformation}\label{ssCob}
The papers \cite{BG} and \cite{Betal} contain a more general class of deformations, with a second parameter of deformation being a function on $\R$. These deformations are obtained by inserting an additional factor into the definition of the map $T_\theta$. We do not need the precise definition of this procedure, see \cite[Section~4.1]{BG} for details, and will only write down the final answer.

Given a smooth function $\tau$ on $\R$ satisfying certain growth conditions, we have a $G$-invariant product $\star_{\theta,\tau}$ on a function space $\E_{\theta,\tau}(G)$ defined by the kernel
$$
K_{\theta,\tau}(x,x')=K_\theta(x,x')\exp\left\{{{\tau\left(\frac{2}{\theta}\sinh(2a)\right) +\tau\left(\frac{2}{\theta}\sinh(-2a')\right)-\tau\left(\frac{2}{\theta}\sinh(2a-2a')\right)}}\right\}.
$$
If $\tau$ is purely imaginary, this kernel defines a unitary $2$-cocycle $\Omega_{\theta,\tau}$ on $\hat G$ such that
$$
\Omega_{\theta,\tau}^*=\int_{G\times G}K_{\theta,\tau}(x,y) \lambda_x\otimes\lambda_y\,dx\,dy\ \ \text{on}\ \ C^\infty_c(G\times G)\subset L^2(G\times G).
$$

In order to understand this cocycle, consider the von Neumann algebra $W^*(\R)$ of $\R$. The conjugation by the inverse of the Fourier transform, defined by
$$
(\F f)(\xi)=\hat f(\xi)=\frac{1}{\sqrt{2\pi}}\int_\R f(t)e^{-i\xi t}dt,
$$
gives an isomorphism $L^\infty(\R)\cong W^*(\R)$, so the unitary $e^{-\tau}\in L^\infty(\R)$ defines a unitary $u_\tau=\F^{-1}e^{-\tau}\F\in W^*(\R)$. Using the embedding $\R\hookrightarrow G$, $t\mapsto(0,0,t)$, we get an embedding $W^*(\R)\hookrightarrow W^*(G)$, so we can consider the unitary $u_\tau$ as an element of $W^*(G)$.

\begin{proposition}
We have $\Omega_{\theta,\tau}=(u_\tau\otimes u_\tau)\Omega_\theta\Dhat(u_\tau)^*$, so the cocycles $\Omega_{\theta,\tau}$ and $\Omega_\theta$ are cohomologous.
\end{proposition}

\bp It is convenient to prove a slightly different statement. Fix functions $f_1,f_2\in S(\R)$. Consider the corresponding elements $b_i=\F^{-1}f_i\F$ of $W^*(\R)\subset W^*(G)$. Explicitly,
$$
b_i=\frac{1}{\sqrt{2\pi}}\int_\R\hat f_i(-t)\lambda_{(0,0,t)}dt.
$$
Consider also the function
$$
\tilde K(x,x')=K_\theta(x,x')f_1\left(\frac{2}{\theta}\sinh(2a)\right)
f_1\left(\frac{2}{\theta}\sinh(-2a')\right)f_2\left(\frac{2}{\theta}\sinh(2a-2a')\right)
$$
and the operator $\tilde\Omega^*=\int_{G\times G}\tilde K(x,y)\lambda_x\otimes\lambda_y\,dx\,dy$, which is at least defined on $C^\infty_c(G\times G)$. We claim that this is a bounded operator on $L^2(G\times G)$ and
\begin{equation}\label{ecob}
\tilde\Omega^*=\Dhat(b_2)\Omega^*_\theta(b_1\otimes b_1).
\end{equation}
Denote $4(\pi\theta)^{-2d-2}A(x,x')$ by $\tilde A(a,a')$ (recall that $A(x,x')$ depends only on the coordinates $a,a'$). Then on $C^\infty_c(G\times G)$ we have:
\begin{align*}
\tilde\Omega^*&=\int\tilde A(a,a')\exp\left\{\frac{2i}{\theta}S(x,x')\right\}
f_1\left(\frac{2}{\theta}\sinh(2a)\right)
f_1\left(\frac{2}{\theta}\sinh(-2a')\right)f_2\left(\frac{2}{\theta}\sinh(2a-2a')\right)\\
&\quad\quad\times\lambda_x\otimes\lambda_{x'}\,dx\,dx'\\
&=\frac{1}{(2\pi)^{3/2}}\int \hat f_1(t_1)\hat f_1(t_1')\hat f_2(t_2)\tilde A(a,a')\exp\left\{\frac{2i}{\theta}\cosh (a)\cosh( a')\omega_0(v ,v')\right\}\\
&\quad\quad\times\exp\left\{\frac{2i}{\theta}\big(t_1\sinh (2a)-t_1'\sinh (2a')+t_2(\sinh (2a)\cosh (2a')-\cosh (2a)\sinh (2a'))\big)\right\}\\
&\quad\quad\times\exp\left\{\frac{2i}{\theta}(\sinh (2a)t'-\sinh (2a')t)\right\}\lambda_x\otimes\lambda_{x'}\,dt_1\,dt_1'\,dt_2\,dx\,dx'\\
&=\frac{1}{(2\pi)^{3/2}}\int \hat f_1(t_1)\hat f_1(t_1')\hat f_2(t_2)\tilde A(a,a')\exp\left\{\frac{2i}{\theta}\cosh( a)\cosh (a')\omega_0(v,v')\right\}\\
&\quad\quad\times\exp\left\{\frac{2i}{\theta}\big(\sinh (2a)t'-\sinh (2a')t\big)\right\}\\
&\quad\quad\times\lambda_{(a ,v ,t-t_1'-t_2\cosh(2a))}\otimes\lambda_{(a' ,v' ,t'-t_1-t_2\cosh(2a'))}dt_1\,dt_1'\,dt_2\,dx\,dx'.
\end{align*}
The group multiplication formula gives
$$
(a ,v ,t-t_1'-t_2\cosh(2a))=
(0,0,-t_2)(a,v,t-t_2\sinh(2a))(0,0,-t'_1).
$$
Hence
\begin{align*}
\tilde\Omega^*
&=\frac{1}{(2\pi)^{3/2}}\int \hat f_1(t_1)\hat f_1(t_1')\hat f_2(t_2)\tilde A(a,a')\exp\left\{\frac{2i}{\theta}\cosh (a)\cosh( a')\omega_0(v ,v')\right\}\\
&\quad\quad\times\exp\left\{\frac{2i}{\theta}\big(\sinh (2a)(t'+t_2\sinh(2a'))-\sinh (2a')(t+t_2\sinh(2a))\big)\right\}\\
&\quad\quad\times(\lambda_{(0,0,-t_2)}\otimes\lambda_{(0,0,-t_2)})(\lambda_x\otimes\lambda_{x'}) (\lambda_{(0,0,-t_1')}\otimes\lambda_{(0,0,-t_1)})dt_1\,dt_1'\,dt_2\,dx\,dx'\\
&=\hat{\Delta}(b_2)\Omega_\theta^*( b_1\otimes b_1).
\end{align*}

Now, choosing a bounded sequence of functions $g_n\in S(\R)$ converging to $e^{\tau}$ pointwise and passing to the limit in identity \eqref{ecob} applied to the pairs $(f_1,f_2)=(g_n,\bar g_n)$, we get the result.
\ep

Therefore, from our perspective, there is no reason to introduce the second deformation parameter~$\tau$, since by Proposition~\ref{pCob} this leads to isomorphic deformations. Note also that on the level of the function spaces~$\E_{\theta,\tau}(G)$ the corresponding $G$-equivariant isomorphism $(\E_\theta(G),\star_\theta)\cong(\E_{\theta,\tau}(G),\star_{\theta,\tau})$ is given by the operator~$\F^{-1}e^{-\tau}\F$.

\bigskip

\section{Open problems} \label{sOP}

The results that we have obtained so far lead to a number of questions that have to be resolved in order to bring the theory to a completely satisfactory level. In this section we list some of the most natural ones.

\subsection{Regularity of cocycles} The problem is to find simple, verifiable conditions for regularity. In view of our considerations in Sections~\ref{ssFixed} and~\ref{ssRegular2} for group duals, it seems unlikely that such conditions exist. At the same time we do not have a single example of a nonregular cocycle on a regular quantum group.

\subsection{Regularity of deformed quantum groups} As has been shown by De Commer~\cite{DC2}, regularity of a quantum group is not preserved under cocycle deformation: the nonregular quantum group~$\tilde E_q(2)$ is obtained by deformation by a cocycle on $SU_q(2)$. In this respect we want to formulate the following question: if $G$ is a regular quantum group and $\Omega$ is a cocycle on $\hat G$, is regularity of $\Omega$ equivalent to regularity of $G_\Omega$? If not, do we have an implication in at least one direction?

\subsection{Generalized fixed point algebras} By Proposition~\ref{pfixed2}, in the case of cocycles on group duals, for the $\Omega$-deformation $A_\Omega$ of a C$^*$-algebra $A$ we have $A_\Omega\subset[\psi_{\beta\otimes\alpha}(\PP^+_{\beta\otimes\alpha})]$. It is not difficult to see that the inclusion can be strict already for $G=\Z$ and $\Omega=1$. What is the proper characterization of elements of $A_\Omega$ in terms of the action $\beta\otimes\alpha$ of $G$ on $C^*_r(\hat G;\Omega)\otimes A$?

Another question is what an analogue of this setting for general quantum groups is. When $G$ is a genuine group, what is of course special, is that any action of $G$ can be viewed as a left or a right action and the tensor product action is always well-defined. When $G$ is a group dual, then again we can always pass from a left to a right action. But in order to define the diagonal action we have to replace the usual tensor product by the braided tensor product $\boxtimes$. As was shown by Yamashita~\cite{Y}, when in addition $G$ is compact, so $C_0(G)=C^*_r(\Gamma)$ for a discrete group $\Gamma$, then $A_\Omega$ is isomorphic to the fixed point algebra $(C^*_r(\Gamma;\Omega)\boxtimes A)^G$. But for general quantum groups, when $(C_0(G),\Delta)$ is neither commutative nor cocommutative, it is not clear to us what the correct analogue of the description of $A_\Omega$ as (a subalgebra of) a fixed point algebra is.

\subsection{Generalization of Rieffel's deformation}
In the setting of Section~\ref{sBG}, Bieliavsky and Gayral defined a $\theta$-deformation $A_\theta$ of any C$^*$-algebra $A$ equipped with an action of $G$. The question is how the algebras $A_\theta$ are related to our algebras $A_{\Omega_\theta}$. As we have seen, this question is not quite trivial already for $A=C_0(G)$. Assuming it can be rigorously settled in this case, for general $A$ both algebras~$A_\theta$ and~$A_{\Omega_\theta}$ can be embedded into $M(C^*_r(\hat G;\Omega_\theta)\otimes A)^{\beta\otimes\alpha}$, and then the question  is whether they coincide. The analogous question for Rieffel's deformation has an affirmative answer~\cite{BNS,N}. There are several reasons why it will be difficult to give a similar proof in the present case. One of them is that in the case of Rieffel's deformation the group~$\R^{2d}$ carrying a dual cocycle was abelian, so the deformation still carried an action of the same group. This is no longer the case for the groups considered by Bieliavsky and Gayral, where we can only hope that $A_\theta$ carries an action of $G_{\Omega_\theta}$.

\bigskip

\end{document}